\documentclass[11pt]{article}

\usepackage{url}
\usepackage{amsfonts,amssymb,amsmath,amsthm}
\usepackage{latexsym}
\usepackage[mathscr]{eucal}
\usepackage{epigraph}

\newtheorem{corollary}{Corollary}[section]

\newtheorem{proposition}{Proposition}[section]

\theoremstyle{definition}
\newtheorem{remark}{Remark}[section]

\numberwithin{equation}{section}

\allowdisplaybreaks[4]
\newcommand{\me}{\mathrm{e}}
\newcommand{\mi}{\mathrm{i}}
\newcommand{\md}{\mathrm{d}}
\begin{document}
\title{Painlev\'{e} Property and Generating Functions for Asymptotics}
\author{A.~V.~Kitaev\thanks{\texttt{E-mail: kitaev@pdmi.ras.ru}}\\
Steklov Mathematical Institute, Fontanka 27,\\ St. Petersburg 191023, Russia}
\date{December 12, 2025}
\maketitle

\epigraphhead[+10pt]{%
\epigraph{\textbf{\textit{To Professor N.\! M.\! Bogoliubov in honour of his 75th birthday}}}{}
}

\begin{abstract}
\noindent
This paper proposes a new approach to the asymptotic analysis of Painlev\'e equations. The approach is based
on representing solutions of the Painlev\'e equations using formal series in two variables,
$\sum_{k=0}^{\infty}y^kA_k(x)$, with rational functions $A_k(x)$. The approach is applied to the asymptotic
analysis of the third degenerate Painlev\'e equation.

\vspace{0.50cm}


\vspace{0.5cm}

\textbf{Abbreviated Title.} Painlev\'e Property and Generating Functions

\vspace{0.5cm}


\end{abstract}
\clearpage
\setcounter{page}{2}
\setlength{\topmargin}{-0.50in}
\section{Introduction} \label{sec:Introduction}
An ordinary differential equation (ODE) is said to have the Painlev\'e property if the positions of singular points,
other than poles, of its general solutions do not depend on the initial values (integration constants) of the
solutions. This definition distinguishes singular points of solutions into movable (poles) and non-movable (branch
points and essential singular) points; in this paper the latter points are often called the regular singular points
and irregular singular points, respectively. When we refer to both types of singular points, regular and
irregular, we call them critical points.

The location of the non-movable singular points can be deduced from the analysis of the form of ODEs; say,
denoting the independent variable of the Painlev\'e equations in the canonical form as $\tau\in\mathbb{C}$,
then it is not complicated to establish that the sixth Painlev\'e equation has three regular (branch)
singularities at $\tau=0$, $\tau=1$, and $\tau=\infty$, the other Painlev\'e equations have one irregular (essential)
singular point at $\tau=\infty$, the fifth and third Painlev\'e equation, additionally, have a regular singular point
at $\tau=0$. The location of the movable poles of general solutions cannot be explicitly found through of any
parameters uniquely define them (initial data, monodromy data), which is a manifestation of their transcendency.
At the same time, the Painlevé property of an ODE imposes certain conditions on the behavior of its solutions in
the neighborhoods of non-movable singular points; this behavior is simpler compared to the behavior of solutions
near singular points of an ODE without the Painlevé property.

The asymptotic behaviour of general solutions of the Painlev\'e equations in the neighbourhood of the non-movable
singular points are also already studied. One of the important consequences of the Painlev\'e property is that it
suggests existence of the connection formulae for the asymptotics of the same solution in different singular
points or between different directions at a given singularity. These connection formulae are already obtained due
to the additional analytic structure which is appeared to be related with the Painlev\'e equations, a theory of
isomonodromy deformations for linear ODEs. It can be said that the leading terms of the asymptotics, together with the
corresponding connection formulae describing the behaviour of the general solution in some neighbourhood of a given
singular point, represent the ``surface'' structure of this
singular point\footnote{\label{foot:surface} This structure can be represented as a neighborhood of a singular point
in $\mathbb{C}$, divided by certain rays into sectors corresponding to the leading terms of the asymptotic series
for general solutions in the sectors. Parametrization of these asymptotic formulae by the monodromy data of the
solutions allows us to describe how the asymptotic behavior varies from sector to sector, which is sometimes called
the nonlinear Stokes phenomenon. Thus, when discussing this structure, we are "looking at the surface."},
so to say, the {\it face} of the singularity.
On the other hand, observing the results concerning the full asymptotic series for the general solutions of the
Painlev\'e equations~\cite{Shimomura1982-2,KimuraH1983,TakanoK1983,TakanoK1986,GILy2013,Shimomura2015}, we conclude
that the Painlev\'e property provides a remarkably simple structure for these series: the reader can verify that
a random addition (or change) of any term of the Painlevé equation that violates the Painlevé property
destroys the structure of these series. This latter structure can be called the ``under the surface''
structure\footnote{\label{foot:subsurface} Clearly, the full asymptotic series contain "too many"
correction terms that cannot be placed on the same surface as the leading term, so we place them under the surface.}
of the singular points, the {\it root} of the singularity, so to speak. The goal of this article is to explicitly
formulate what we mean by the root structure, develop an analytical tool for studying it, and illustrate its main
features using the example of the degenerate third Painlev\'e equation.

Let $u(\tau)$ be a solution of an equation with the Painlev\'e property and $\tau_0$ its critical point,
denote $\Omega\equiv\Omega_{\tau_0}$ a connected domain in the complex $\tau$-plane such that $\tau_0$  belongs
to its boundary, $\tau_0\in\partial\Omega$. We call $\Omega$ a proper domain at $\tau_0$ if
(1) $\tau_0$ is a regular singular point and $\Omega$ is an open disc centered at $\tau_0$, with a radius cut,
if the function $u(\tau)$ has a countable number of poles/zeros accumulating at $\tau_0$, then additionally, we
have to delete from $\Omega$ closed discs centred at these poles/zeros with an appropriately shrinking with
$\tau\to\tau_0$ radiuses (see, e.g. a construction presented in \cite{KitVar2025}), so that in the presence of such
poles/zeros $\Omega$ is a multiply connected domain,
(2) $\tau_0$ is an essential singularity, $\Omega$ is a``fat'' ray, i.e., the ray with the origin at $\tau_0$
together with some shrinking as $\tau\to\tau_0$ neighbourhood. In case there are sequences of poles/zeros accumulating
at $\tau=\tau_0$ applies a comment similar to the one in item (1). The exact definition of $\Omega$ depends on
a particular equation with the Painlev\'e property under investigation and $\arg\tau$ and will be specified below
for the solutions of the degenerate Painlev\'e equation. All examples that are considered by the author so far
are related with the Painlev\'e equations in the canonical form. For these equations a regular singular point, if any,
is always located at $\tau_0=0$ and the essential singular point, if it exists, at $\tau_0=\infty$. As mentioned above,
the only Painlev\'e equation that does not possess an irregular singularity is the sixth Painlev\'e equation,
instead it has at $\tau=\infty$ a regular singular point and one additional regular singularity at $\tau=1$.

Now, the root structure of the non-movable critical points can be defined in the following way.
Let $u(\tau)$ be any solution of some Painlev\'e equation and $\tau_0$ its critical point,
then there exist a proper domain $\Omega_{\tau_0}$, functions $x=x(\tau)$ and $y=y(\tau)$ holomorphic in
$\Omega_{\tau_0}$ and independent of the initial data of $u(\tau)$, such that
\begin{equation}\label{eq:uAxy}
u(\tau)\underset{\genfrac{}{}{0pt}{}{\tau\to\tau_0}{\tau\in\Omega_{\tau_0}}}{=}\tau^rA(x(\tau),y(\tau)),
\end{equation}
where $r$ is a rational number which depends on the Painlev\'e equation and $\tau_0$, and
$A(x,y)$ is a formal expansion,
\begin{equation}\label{eq:A-def}
A(x,y)=\sum_{k=0}^{\infty}y^kA_{k}(x),
\end{equation}
where $A_k(x)$ are rational functions of $x$. The formal expansion~\eqref{eq:A-def} after substituting $x=x(\tau)$
and $y=y(\tau)$ becomes an asymptotic series in $\tau$, possibly after specifying the parameters they may have.

This structure should be supplemented with the following, less obvious features: (1) the choice of the functions
$x(\tau)$ and $y(\tau)$ is not unique; (2) all asymptotic series possess a certain symmetry which can be
exploited by making the corresponding symmetry transformations for the functions $x(\tau)\to \tilde{x}(\tau)$ and
$y(\tau)\to\tilde{y}(\tau)$, and substituting the tilde-functions for $x(\tau)$ and $y(\tau)$ in the
equation~\eqref{eq:A-def}. We get also a valid asymptotic series, which can be, in particular, used together
with the original series~\eqref{eq:A-def} to obtain a uniform asymptotic expansion; and (3) there is always
a proper choice of the functions  $x=x(\tau)$ and $y=y(\tau)$ such that the following representation is valid
\begin{equation}\label{eq:AxyB}
A(x,y)=\sum_{k=0}^{\infty}x^kB_{k}(y),
\end{equation}
with some rational functions $B_k(y)$; We call the varables $x(\tau)$ and $y(\tau)$  and the
expansions~\eqref{eq:A-def} and \eqref{eq:AxyB} {\it conjugate} functions and, respectively, expansions.
(4) The poles of the functions $A_k(x)$ and
$B_k(y)$, $k\in\mathbb{N}$ are defined by the functions $A_0(x)$ and $B_0(y)$, respectively.

We claim that this root structure holds for all asymptotics at all critical points of all Painlevé equations.
The functions $x(\tau)$ and $y(\tau)$, as well as the formal expansion $A(x,y)$ depend on each object of study
mentioned in the previous sentence.

The reader noticed that the functions $x(\tau)$ and $y(\tau)$ appear in a pair, so that they can be interchanged
and, thus, it is not clear whether there is any difference between expansions~\eqref{eq:uAxy} ($A-expansion$) and
\eqref{eq:AxyB} ($B-expansion$)? In fact, we distinguish these expansions by referencing an original asymptotic
expansion, which is an external object for us, which was obtained by other authors as $B$-expansion (base expansion).
In all main examples for the general solutions, that we studied so far, the $B$-expansion is symmetric,
$y^kB_k(x)=\tilde{y}^kB_{k}(\tilde{x})$ and $B_k(y)$ are polynomials. These $B$-expansions serves us the sources
of the functions $x(\tau)$ and $y(\tau)$. These features can be readily observed from Section~\ref{sec:trig}.
After the proper variables $x(\tau)$ and $y(\tau)$ are defined one can obtain, with a help of the Painlev\'e
equation under investigation, a partial differential equation (PDE) for determining $A/B$-expansions.

In principle, a reader familiar with Appendices A-C \cite{KitVar2025} and who has studied Sections~\ref{sec:log}
and \ref{sec:trig} of this paper should have no problem finding the variables $x(\tau)$ and $y(\tau)$ for the
``power'' and/or ``logarithmic'' asymptotics at regular singularities and the ``trigonometric'' and/or
``trancated'' asymptotics at irregular singular points for any Painlev\'e equation, provided that the leading term
of the corresponding asymptotics is known. The situation with ``elliptic'' asymptotics at irregular singularities is
more complicated because we have no one example of the complete asymptotic expansion for this type of the asymptotic
behaviour. Two sighting shots for elliptic expansions are made in Section~\ref{sec:elliptic}, both were partially
successful, because in Subsection~\ref{subsec:e-non-conjugate}, we introduce variables $x(\tau)$ and $y(\tau)$ such
that the expansion~\eqref{eq:A-def} is valid but the variables $x(\tau)$ and $y(\tau)$ are non-conjugate, in
Subsection~\eqref{subsec:e-conjugate} we introduced formally conjugate variables $x(\tau)$ and $y(\tau)$ and not
to overload the notation we denote one of the expansions as $B$-expansion but in this case the conjugate variables
appear to be symmetric, i.e., $\tilde{x}(\tau)=y(\tau)$ and $\tilde{y}(\tau)=x(\tau)$, so that in this case,
these conjugate variables define a symmetric expansion, rather than the conjugate one. The main propery of
the $B$-expansion is violated, it is not symmetric. An inderect consequence of it is the presence an infinite
sequence of constants of integration in the expansion, while a proper $B$-expansion should not have any indetermined
coefficients and we use it for determination of initial conditions for
the functions $A_k(x)$.\footnote{\label{foot:B} Calculating the $B$-expansion at first glance seems similar to
substituting the asymptotic ansatz into the Painlevé equation. However, from a computational standpoint,
there is a significant difference. When substituting the asymptotic ansatz into the Painlev\'e equation, we must
initially allocate memory for all its terms with undetermined coefficients and perform the calculations by storing
all these terms, as well as the corresponding algebraic equations, in memory, then solving them sequentially.
Meanwhile, the $B_k(y)$ function appears in explicit (complex) form only after its calculation, and the memory usage
increases step by step. Furthermore, calculating the coefficients hidden in $B_k(y)$ is done by integration, while
calculating the same set of coefficients in the asymptotic ansatz is done by differentiation. For asymptotic series
with a large number of terms, the advantage of using the $B$-expansion becomes obvious.}
In theory, we could probably determine these constants by substituting this expansion into the
Painlev\'e equation, but that is beyond the scope of our method: leave the Painlevé equation itself alone.

Let us recall that modern methodologies based on isomonodromy deformations do not require a preliminary asymptotic
analysis of the Painlev\'e equations and provide the researcher with the "face structure" of singularities.
The PDEs for A/B expansions, which determine the ``root structure'' of singularities, are, of course, based on
the corresponding Painlevé equations, but ideologically, the method is different. While in standard classical
asymptotic analysis we study a simplified/truncated/reduced version of the Painlevé equation, which defines the
leading term of the asymptotics; for general solutions, this is the equation for elliptic functions or their
degeneracies, in our approach we deal only with rational solutions and avoid any asymptotic ansatzes and expansions
that appeared in studies based on the study of simplified equations.

In this paper, we discuss the application of A/B expansions to the study of the asymptotic behavior of the Painlev\'e
equations, which has a dual nature:
(1) we show that the rational functions $A_k(x)$ and $B_k(y)$, $k\in\mathbb{Z}_{\geqslant0}$ generate explicit
formulae for infinite/finite subsequences of coefficients of the asymptotic expansions therefore we call them
generating functions. The formal expansions~\eqref{eq:A-def} and \eqref{eq:AxyB}, representing collections of
these generating functions, we call {\it supergenerating functions};
(2) the $A$- and $B$-expansions for a suitable  choice of the parameters in the functions $x=x(\tau)$ and
$y=y(\tau)$ define asymptotic expansions, which are an analytic continuation of each other. Moreover, these
expansions have an overlapping domains of validity in the parameter space. This fact has the following consequence,
if a parametrization of one of the asymptotic expansions of the solution $u(\tau)$ by its monodromy data is obtained,
then the same parametrization is true for the other asymptotic expansion, so that mo any additional calculations
for finding of such monodromy parametrizations are required for the latter asymptotics (see example
in Section~\ref{sec:trig}).

In this paper $u(\tau)$ is a solution of the degenerate third Painlev\'{e} equation
\begin{equation} \label{eq:dp3u}
u^{\prime \prime}(\tau) \! = \! \frac{(u^{\prime}(\tau))^{2}}{u(\tau)} \! - \! \frac{u^{\prime}(\tau)}{\tau}
\! + \! \frac{1}{\tau} \! \left(-8 \varepsilon (u(\tau))^{2} \! + \! 2ab \right) \! + \! \frac{b^{2}}{u(\tau)},
\end{equation}
where the prime denotes differentiation with respect to $\tau$, $\varepsilon=\pm1$,
$b\in\mathbb{R}\setminus\lbrace 0\rbrace$ is a scaling parameter, and $a\in\mathbb{C}$ is a so-called formal
monodromy parameter, which cannot be changed by any scaling transformations. Equation~\eqref{eq:dp3u} is also
referred to as the third Painlev\'e equation of type $D_7$ (see \cite{OKSO2006}). The choice of this equation 
is due to
the fact that it has both types of critical points, has one (non-scaling) parameter $a$, but not as many such
parameters as other Painlev\'e equations with two types of singular points, and, in fact, it is somewhat more
involved as an object of asymptotic study than other Painlev\'e equations, which look even more complicated.

Finally, let us discuss some further questions and perspectives related with the root structure of the
critical singularities.
(1) The representation~\eqref{eq:uAxy} for asymptotcs of solutions is an implication of the Painlev\'e property
for the second order ODEs, i.e., necessary conditions for the Painlev\'e property of the second order ODEs.
So, a natural question is what is a generalization of these expansions for the higher order Painlev\'e equations
of one and several variables?

(2) Is that possible that a moving critical singularity of general solution of a second order ODE has a root structure
described by the expansion~\eqref{eq:uAxy}?

(3) If some second order ODE have some non-movable critical points with the root structure~\eqref{eq:uAxy} and
some movable critical points, then the movable critical points cannot enter into some neighbourhoods of the
non-movable critical points because they would destroy (it contrudicts to) the structure~~\eqref{eq:uAxy}.
So, the union of neighbourhoods of the non-movable critical points with the root structure~\eqref{eq:uAxy}
form a ``repulsive'' domain for the movable critical points, i.e., a domain without movable critical points,
this domain can be called the Painlev\'e domain and its completion is the non-Painlev\'e domain where the
movable critical points are walking. So, this idea gives a theoretical explanation of a possible existence of
the equations with partial Painlev\'e property. One of the candidates to equations of such type were considered
in \cite{Kit2012}. The major question here is that although the root structure~\eqref{eq:uAxy} is a local property
of a critical point, it is defined by all terms of the corresponding ODE. Is that possible that such ODE which
holds the structure of the non-movable points allows to move some critical singular points?

\section{Regular Singularity: Logarithmic Asymptotics}\label{sec:log}
Equation~\eqref{eq:dp3u} has one regular singular point located at $\tau=0$. In the paper~\cite{KitVar2025} we studied
the one-parameter solution of this equation which has the following asymptotic
expansion\footnote{\label{foot:ln-reg-sol} For details, see Appendix C of \cite{KitVar2025}.
In this section, for the convenience of the reader, we retain the notation adopted in \cite{KitVar2025},
in particular, in this section $\varepsilon=+1$, to obtain formulae with $\varepsilon=+/-1$ it is sufficient to
replace $u->\varepsilon u$ and $b->\varepsilon b$ in both parts of all equations. In the equations with
the function $A(x,y)$, the function $A(x,y)$ should be changed to $\varepsilon A(x,y)$.}
\begin{equation}\label{app:eq:0-u-expansionLog2}
\begin{gathered}
u(\tau)=\sum_{k=0}^{+\infty}\tau^{2k-1}\sum_{m=-2\lfloor k/2\rfloor}^{+\infty}\tilde{c}_{2k-1,m}(\ln\tau)^{-m},\\
\tilde{c}_{-1,0}=\tilde{c}_{-1,1}=0,\quad
\tilde{c}_{-1,2}=-\frac14.
\end{gathered}
\end{equation}
The expansion~\eqref{app:eq:0-u-expansionLog2} depends on the single parameter $\tilde{c}_{-1,3}\in\mathbb{C}$ and
is convergent in a neighbourhood of $\tau=0$; furthermore, as in \cite{KitVar2025}, we assume that
$|\arg\,\tau|<\pi$ and the principle branch of $\ln$-function is chosen.
The expansion~\eqref{app:eq:0-u-expansionLog2} is valid for all values of the parameter of formal monodromy
$a\in\mathbb{C}$.
We say that the coefficients $\tilde{c}_{2k-1,m}$ are the \emph{coefficients of level $k$}. In contrast with
the other types of the asymptotic expansions as $\tau\to0$ of solutions $u(\tau)$ the
number of coefficients that belong to each level for $\tilde{c}_{-1,3}\neq0$ is infinite. In Appendix C of
\cite{KitVar2025}, we constructed the super-generating function $\tilde{A}(x,y)$ with $x=1/\ln(\tau)$ and $y=\tau^2$,
that allowed us to find explicit formulae for the coefficients $\tilde{c}_{2k-1,m}$ for all
$m\geqslant-2\lfloor k/2\rfloor$ of the levels $k=0,1,2$, and $3$.
Using the super-generating function $\tilde{A}(x,y)$, one can continue to construct the coefficients
$\tilde{c}_{2k-1,m}$ for higher levels $k$, however the corresponding formulae quickly become rather complicated,
nevertheless, some general properties of the coefficients can be derived. By other words the construction presented
in Appendix C.2 says that for any given $k$ the corresponding inner series of the double
series~\eqref{app:eq:0-u-expansionLog2} is nothing but a rational generating function $\tilde{A}_k(x)$.

The choice of the functions $x(\tau)=1/\ln(\tau)$ and $y(\tau)=\tau^2$ is not a proper choice in the sense discussed
in the Introduction. Because, the conjugate expansion with respect to variable $x$ is a Larent (doubly infinite)
expansion, and does not have the form~\eqref{eq:AxyB}. Below we consider the proper choices of the pair of functions
$x(\tau)$ and $y(\tau)$.

The first choice is
\begin{equation}\label{eq:xy-newdef-regular}
x(\tau)=\tau^2\ln(\tau),\qquad
y(\tau)=1/\ln(\tau).
\end{equation}
Using variables~\eqref{eq:xy-newdef-regular} the expansion~\eqref{app:eq:0-u-expansionLog2} can be
written in terms of the function $A(x,y):=\tau u(\tau)$ as the expansion~\eqref{eq:uAxy}, where
\begin{equation}\label{eq:Ak-def-regular}
A_k(x)=\sum_{n=0}^{\infty} \tilde{c}_{2n-1,k-n}x^n, \qquad
k\geqslant0.
\end{equation}
In equation~\eqref{eq:Ak-def-regular} we assume that for $k=0$ all coefficients corresponding to odd
$n\in\mathbb{N}$ vanish, i.e., $\tilde{c}_{4l-3,1-2l}=0$ for $l\in\mathbb{N}$.
Define the linear differential operator $D$ acting in the space of formal power
series of two variables $x$ and $y$ as follows:
\begin{equation}\label{app:eq:tildeD-definition}
D=(y+2)x\frac{\partial}{\partial x}-y^2\frac{\partial}{\partial y};
\end{equation}
then, the function $A\equiv A(x,y)$ solves the PDE
\begin{equation}\label{app:eq:tildeA-PDE}
D^2(\ln A)=-8A+2a\frac{bxy}{A}+\left(\frac{bxy}{A}\right)^2.
\end{equation}
Now we substitute the expansion~\eqref{eq:uAxy} for the function $A$ into the equation~\eqref{app:eq:tildeA-PDE},
equate the coefficients $y^k$ to zero and obtain differential equations for determining the generating functions
$A_k(x)$. For $A_0(x)$ we obtain a nonlinear differential equation with the only rational solution
\begin{equation}\label{eq:A0}
A_0(x)=\frac{Cx^2}{(1 + Cx^2/4)^2},
\end{equation}
Expanding $A_0(x)$ into the Taylor series, and comparing the result with the expansion~\eqref{eq:Ak-def-regular},
one finds that $C=\tilde{c}_{3,-2}$ and for $l\in\mathbb{N}$ finds:
\begin{align}
&n=2l-1:&
&\tilde{c}_{4l-3,1-2l}=0,\qquad
l\in\mathbb{N},&
\label{eqs:c4l-3,1-2l}\\
&n=2l:&
&\tilde{c}_{4l-1,-2l}=\left(-\frac{1}{4}\right)^{l-1} l\tilde{c}_{3,-2}^l
=-\frac{lb^{2l}(a^2+1)^l}{4^{2l-1}}.&\label{eqs:c4l-1,-2l}
\end{align}
In the last equation of \eqref{eqs:c4l-1,-2l}, we used the value of $\tilde{c}_{3,-2}$ which was calculated
in \cite{KitVar2025} via a direct substitution of the expansion~\eqref{app:eq:0-u-expansionLog2} into the
equation~\eqref{eq:dp3u}. Other generating functions $A_k(x)$, $k\in\mathbb{N}$, are defined as rational solutions
of non-homogeneous linear equations of the second order, the homogeneous parts of which are degenerate hypergeometric
equations. In this case appear an additional problem, that does not occur in other situations that we considered
so far, namely, this is a resaunant case, by other words, the homogeneous part of the ODEs defining $A_k(x)$ have
a rational solution, so that the rational solution that define $A_k(x)$ is not unique and we faced to the
problem of determining for each $k$ one initial parameter. This parameter, of course, can be calculated but it
requires additional efforts. Let's see what happens for $A_1(x)$;
\begin{equation}\label{eq:A1CC1}
\begin{aligned}
A_1(x)&=\frac{C_1x^2(Cx^2-4)}{(Cx^2+4)^3}\\
&+\frac{x\big(ab(Cx^2-12)(C^2x^4+56Cx^2-48)+2304Cx\big)}{18(Cx^2+4)^3}
\end{aligned}
\end{equation}
The first term in equation~\eqref{eq:A1CC1} is a rational solution of the homogeneous ODE, depending on the constant
of integration $C_1$, which must be determined, and the parameter $C$ is the same as in equation~\eqref{eq:A0},
which is already known. Expanding $A_1(x)$ into the Taylor series one obtains
\begin{equation}\label{eq:A1Taylor}
A_1(x)=\frac{ab}2x+\left(2C-\frac{C_1}{16}\right)x^2-abCx^3+\mathcal{O}\big(x^4\big).
\end{equation}
Comparing the expansion~\eqref{eq:Ak-def-regular} for $k=1$ with the expansion~\eqref{eq:A1Taylor} for $n=0,1,2$,
and $3$, one finds
\begin{equation}\label{eq:coeffsA1}
\tilde{c}_{-1,1}=0,\;\;
\tilde{c}_{1,0}=\frac{ab}{2},\;\;
\tilde{c}_{3,-1}=-\frac{C_1}{16}+2C,\;\;
\tilde{c}_{5,-2}=-abC=\frac{ab^3(a^2+1)}{4},
\end{equation}
The coefficient $\tilde{c}_{3,-1}$ was calculated in \cite{KitVar2025},
\begin{equation}\label{eq:c3m1}
\tilde{c}_{3,-1}=-b^2\big((a^2+1)\tilde{c}_{-1,3}-a^2-1/2\big)=4C(\tilde{c}_{-1,3}-1)-b^2/2.
\end{equation}
We continue to use the notation $C=-b^2(a^2+1)/4$, since in terms of $C$ the formulae look simpler.
Substituting instead of $\tilde{c}_{3,-1}$ , its expressions, given by the equations~\eqref{eq:c3m1},
into the third equation of the list~\eqref{eq:coeffsA1}, one gets
\begin{equation}\label{eq:C1}
C_1=8b^2\big(2(a^2+1)\tilde{c}_{-1, 3}-3a^2-2\big)=-8(8C\tilde{c}_{-1,3}-12C-b^2).
\end{equation}
Now, we can decompose $A_1(x)$ into partial fractions,
\begin{equation}\label{eq:A1parfrac}
A_1(x)=\frac{ab}{18}+\sum_{k=1}^3\frac{\xi_{k,1}x+\xi_{k,0}}{(1+Cx^2/4)^k},
\end{equation}
where
\begin{gather*}
\xi_{1,1}=\frac{4}{9}ab,\quad
\xi_{2,1}=-\frac{32}{9}ab,\quad
\xi_{3,1}=\frac{32}{9}ab,\quad\label{eq:ksi11-ksi31}\\
\xi_{1,0}=-16\tilde{c}_{-1,3}+\frac{2b^2}{C}+24,\quad
\xi_{2,0}=48\tilde{c}_{-1,3}-\frac{6b^2}{C}-64,\\
\xi_{3,0}=-32\tilde{c}_{-1,3}+\frac{4b^2}{C}+40.\label{eq:ksi10-ksi30}
\end{gather*}
Developing $A_1(x)$ given in terms of the partial fractions~\eqref{eq:A1parfrac} into the Taylor series at  $x=0$
and comparing the coefficients of the latter expansion with the coefficients of expansion~\eqref{eq:Ak-def-regular}
for $k=1$, one arrives at the following explicit formulae,
$\tilde{c}_{1,0}=ab/18+\xi_{1,1}+\xi_{2,1}+\xi_{3,1}=ab/2$ and for $l\in\mathbb{N}$:
\begin{align}
&n=2l+1:\;
\tilde{c}_{4l+1,-2l}=\left(\frac{-C}{4}\right)^l
\left(\xi_{1,1}+(l+1)\xi_{2,1}+\frac{(l+1)(l+2)}{2}\xi_{3,1}\right),
\label{eqs:c4l+1,-2l}\\
&n=2l:\;\;\;\;\;\,
\tilde{c}_{4l-1,1-2l}=\left(\frac{-C}{4}\right)^l
\left(\xi_{1,0}+(l+1)\xi_{2,0}+\frac{(l+1)(l+2)}{2}\xi_{3,0}\right).
\label{eqs:c4l-1,1-2l}
\end{align}
Note that formula~\eqref{eqs:c4l-1,1-2l} works also for $l=0$, $\tilde{c}_{-1,1}=\xi_{1,0}+\xi_{2,0}+\xi_{3,0}=0$.
\begin{remark}
The determination of the functions $A_k(x)$ for $k\geqslant2$ is similar, but the corresponding expressions quickly
become quite complex.
I calculated $A_2$, $A_3$, and $A_4$ using \textsc{Maple} and did not notice a significant increase in computation
time, so I suspect that in practice, several more generating functions could be calculated, but \textsc{Maple}
might not be able to display them.
To determine the integration constants of type $C_1$ (see above) for functions $A_k(x)$, it is necessary to know
the corresponding coefficient of level $3$, $\tilde{c}_{3,k-2}$; these coefficients (for all $k$) are calculated
in Appendix C.2 \cite{KitVar2025} (see equation~(C.30) in \cite{KitVar2025}).

The order of pole of the function $A_k(x)$ is $k+2$,
so that the corresponding coefficients $\tilde{c}_{2n-1,k-n}=\mathcal{O}(n^{k+1})$.
\hfill$\blacksquare$\end{remark}

Consider the conjugate expansion~\eqref{eq:AxyB}. The corresponding generating functions $B_k(y)$ read,
\begin{equation}\label{eq:Bkreg}
B_k(y)=\sum_{n=0}^{\infty} \tilde{c}_{2k-1,n-k}x^n, \qquad
k\geqslant0.
\end{equation}
Comparing \eqref{eq:Bkreg} with the equation (C.12) in Appendix C.2 of \cite{KitVar2025}, we see that the set of
coefficients defining $B_k$ is the same as the one for the generating function $\tilde{A}_k$ considered in
\cite{KitVar2025}, so that the construction below should produce the same explicit expressions for the coefficients
$\tilde{c}_{2k-1,m}$ as those obtained in Appendix C.2 of \cite{KitVar2025}. This construction, however, is
important, because it illustrates the idea of choosing a proper pair of functions $x(\tau)$ and $y(\tau)$
and the existence of the conjugate expansions discussed in the Introduction. Since in this case the construction of
the expansion~\eqref{eq:AxyB} is similar that considered in \cite{KitVar2025}, we will only briefly outline its
main features below.
\begin{equation}\label{eq:B0reg}
B_0(y)=-\frac{1}{\left(1/y-C\right)^2},\qquad
C=-2\tilde{c}_{-1,3}.
\end{equation}
Function $B_0(y)$ with the change $y\to x$ coincide with function $\tilde{A}_0(x)$ introduced in Appendix C.2 of
\cite{KitVar2025}. The functions $B_k(y)$ for $k\geqslant$ are determined via substitution of the
expansion~\eqref{eq:AxyB} into the equation~\eqref{app:eq:tildeA-PDE} and equating successively the coefficients
of $y^k$ to zero. Like in the case with functions $A_k(x)$, considered above, one obtains nonhomogeneous ODEs
for the determinantion of the functions $A_k(x)$ whose inhomogeneous part is a rational function of $y$, and the
homogeneous part is the second order ODE equivqlent to a degenerate hypergeometric equation. In contrast with
the case considered above the homogeneous  ODE does not have rational solutions, and the determining of the
functions $B_k(y)$ is straightforward: it is the only rational solution of the inhomogeneous ODE.
The final answer is $B_k(y)=y^k\tilde{A}_k(x)|_{x=y}$, where $\tilde{A}_k(x)$ are the functions determined in
Appendix C.2 of \cite{KitVar2025}; this relation means that the conjugate $B$-expansion coincide with
the expansion studied in Appendix C.2 of \cite{KitVar2025}.

All the expansions of \eqref{app:eq:0-u-expansionLog2}, $A$, and $B$ considered in this section are asymptotic
series. From a practical point of view, the series ~\eqref{app:eq:0-u-expansionLog2} is difficult to use; the
expansion of $A$ converges very slowly; the expansion of $B$ is the best of them, but the formulas for the functions
$B_k(y)$ $k>6$ become very complicated (see the discussion at the end of Appendix C.2 in \cite{KitVar2025}), so
further research is required even to display the results of their calculations.

\section{Irregular Singularity: Trigonometric Asymptotics}\label{sec:trig}
Consider asymptotic expansion for the general solution of equation~\eqref{eq:dp3u},
\begin{equation}\label{eq:app:u-asympt-regular-complete}
\begin{gathered}
u(\tau)\underset{\genfrac{}{}{0pt}{}{\mathrm{Re}\,\tau\to+\infty}{\tau\in\Omega}}{=}
\frac{\varepsilon(\varepsilon b)^{2/3}}{2}\tau^{1/3}\left(1+
\sum_{k=1}^\infty\frac{1}{\tau^{k/3}}\sum_{j=-k}^{j=k}a_{k,j}w^{j}\right),\\
\Omega:=\{|\mathrm{Im}\,\theta(\tau)|<\delta,\;\delta>0\},
\end{gathered}
\end{equation}
where
\begin{gather}
w:={\tau}^{\frac23(\tilde{\nu}+1)}\me^{\mi\theta(\tau)},\quad
\theta(\tau):=3^{3/2}(\varepsilon b)^{1/3}\tau^{2/3},\label{eq:infty:regular-w-theta-tau}\\
|\mathrm{Re}\,\tilde{\nu}+1|<1/2,\label{eq:Omega-infty-def}
\end{gather}
and
\begin{equation}\label{eq:coeffs:a1j}
a_{1,0}=0,\qquad
a_{1,1}\,a_{1,-1}=-\frac{\mi(\tilde{\nu}+1)}{\sqrt{3}(\varepsilon b)^{1/3}}.
\end{equation}
The remaining coefficients $a_{k,j}$ can be determined uniquely in terms of $a_{1,1}$, $a_{1,-1}$, and the
coefficients of the equation~\eqref{eq:dp3u} by substituting the expansion~\eqref{eq:app:u-asympt-regular-complete}
into it.
In Appendix C \cite{KitVar2023} we discussed various properties of this expansion, in particular, we
conjectured explicit formulae for $a_{k,\pm k}$ and $a_{k,\pm(k-2)}$, and we also mentioned that these formulae
can be proved using generating function technique, but the details were omitted. In this paper, we provide
the corresponding construction. Here, we follow the notation introduced in \cite{KitVar2023},
\begin{equation}\label{eqs:varkappa-alpha}
\varkappa=\tilde{\nu}+1,\qquad
\alpha=2\mi\sqrt{3}\,a.
\end{equation}
Denote,
\begin{equation}\label{eq:xy-infty1}
x(\tau)=\theta(\tau)^{-1/2}w^{-1},\quad
y(\tau)=w,\qquad
\theta(\tau)^{1/2}=3^{3/4}(\varepsilon b)^{1/6}\tau^{1/3}=:\Theta\tau^{1/3},
\end{equation}
where the function $w$ is defined in \eqref{eq:app:u-asympt-regular-complete}.
\subsection{Super-generating function $\overset{\;\infty}{A}(x,y)$}\label{subsec:Ainfty:trig}

The super-generating function $\overset{\;\infty}{A}(x,y)$ is a function of two
independent variables $x$ and $y$,
\begin{equation}\label{eq:A-inftyDef}
\overset{\;\infty}{A}(x,y)=\sum_{k=0}^{\infty}y^k\overset{\;\infty}{A}_k(x),
\end{equation}
where the generating functions
\begin{equation}\label{eq:Ak-inftyDef}
\overset{\;\infty}{A}_k(x)=\sum_{n=\left\lfloor\tfrac{k+1}{2}\right\rfloor}^{\infty}\Theta^na_{n,k-n}x^n,\qquad
k\in\mathbb{Z}_{\geqslant0},\qquad
a_{0,0}=1.
\end{equation}
Assuming that the function $\stackrel{\;\infty}{A}(x,y)$ is constructed, the asymptotics of the solution $u(\tau)$
can be written in the form
\begin{equation}\label{eq:u-Ainfty}
u(\tau)=\frac{\varepsilon(\varepsilon b)^{2/3}}{2}\tau^{1/3}\stackrel{\;\infty}{A}(x(\tau),y(\tau)),
\end{equation}
where $x(\tau)$ and $y(\tau)$ are defined in \eqref{eq:xy-infty1}.
The linear differential operator $D$ acting in the space of formal power series in two variables $x$ and $y$
in this case is defined as
\begin{equation}\label{eq:D-inftyDef}
D=-\frac{x}{2}\frac{\partial}{\partial x}+\left(\varkappa+\frac{\mi}{x^2y^2}\right)
\left(y\frac{\partial}{\partial y}-x\frac{\partial}{\partial x}\right).
\end{equation}
The function $\overset{\;\infty}{A}\equiv\overset{\;\infty}{A}(x,y)$ satisfies the following PDE
\begin{equation}\label{eq:A-inftyPDE}
D^2\ln\,\overset{\;\infty}{A}=-\frac{1}{3x^4y^4}\left(\overset{\;\infty}{A}
-\frac{1}{\Big(\overset{\;\infty}{A}\Big)^2}\right)-\frac{\mi\alpha}{2x^2y^2\overset{\;\infty}{A}}.
\end{equation}
Substituting the expansion~\eqref{eq:A-inftyDef} into the equation~\eqref{eq:A-inftyPDE} and equating to zero the
coefficients at the powers $y^k$, $k=-4,-3,\ldots$, one obtains for $k=-4$ a nonlinear ODE for the generating
function $A_0(x)$,
\begin{equation}\label{eq:A0diff}
\left(x\frac{\md}{\md x}\right)^2\ln\overset{\;\infty}{A}_0(x)=
\frac13\left(\overset{\;\infty}{A}_0(x)-\frac{1}{\Big(\overset{\;\infty}{A}_0(x)\Big)^2}\right),
\end{equation}
and linear second order ODEs for the functions $\overset{\;\infty}{A}_k(x)$, $k\geqslant1$. The general solution of
equation~\eqref{eq:A0diff} can be written using the Weierstass \text{\Large{$\wp$}}-function of $\ln x$,
However, according to what is written in the Introduction, there must be a rational solution to this equation
and it does exist,
\begin{equation}\label{eq:A0-inftyRat}
\overset{\;\infty}{A}_0(x)=1+\frac{b_{1,-1}x}{\left(1-b_{1,-1}x/6\right)^2},
\end{equation}
where $b_{1,-1}\in\mathbb{C}$ is the constant of integration. By expanding the function $\overset{\;\infty}{A}_0$,
given by the equation~\eqref{eq:A0-inftyRat}, into a Taylor series at x=0 and comparing the result with
expansion~\eqref{eq:Ak-inftyDef} for $k=0$, we obtain formulae for the coefficients $a_{n,-n}$,
\begin{equation}\label{eq:ak-k}
a_{n,-n}=\frac{n(b_{1,-1}\Theta^{-1})^n}{6^{n-1}}=\frac{na_{1,-1}^n}{6^{n-1}}.
\end{equation}
The function $\overset{\;\infty}{A}_1(x)$ solves the second order linear homogeneous ODE,
\begin{equation}\label{eq:A1diff}
\begin{aligned}
&x(xb_{1,-1}-6)(x^2b_{1,-1}^2+24xb_{1,-1}+36)\frac{\md^2\overset{\;\infty}{A}_1(x)}{\md x^2}\\
&=(x^3b_{1,-1}^3-54x^2b_{1,-1}^2
-540xb_{1,-1}-216)\frac{\md\overset{\;\infty}{A}_1(x)}{\md x}\\
&+72b_{1,-1}(xb_{1,-1}+6)\overset{\;\infty}{A}_1(x).
\end{aligned}
\end{equation}
One of the independent solutions of \eqref{eq:A1diff} is a rational function and the other includes $\ln(x)$,
so that
\begin{equation}\label{eq:A1-rational}
\overset{\;\infty}{A}_1(x)=\frac{C_1x^2(xb_{1,-1}+6)}{(xb_{1,-1}-6)^3}=-\frac{C_1}{36}x^2+\mathcal{O}(x^3),
\end{equation}
where $C_1$ is the integration constant. Comparing the expansion~\eqref{eq:A1-rational} with the
expansion~\eqref{eq:Ak-inftyDef} for $k=1$, we confirm that $a_{1,0}=0$, and find the relation
$-C_1/36=a_{2,-1}\Theta^2$. In \cite{KitVar2023} we found that $a_{2,\pm1}=0$, therefore
$C_1=0$, this, in turn, means that $a_{n,1-n}=0$ for all $n\in\mathbb{N}$.
The function $\overset{\;\infty}{A}_2(x)$ is a rational solution of a linear second order inhomogeneous ODE.
This solution can be written as
\begin{equation}\label{eq:A2-C2}
\begin{gathered}
\overset{\;\infty}{A}_2(x)=\frac{C_2x^3(xb_{1,-1}+6)}{(xb_{1,-1}-6)^3}-\frac{\mi x}{12b_{1,-1}(xb_{1,-1}-6)^3}
((\varkappa + 1)b_{1,-1}^5x^5\\
+ 6(\alpha + 9\varkappa + 9)b_{1,-1}^4x^4
- 108(4\alpha + 53\varkappa + 53)b_{1,-1}^3x^3\\
- 216(6\alpha + 143\varkappa + 151)b_{1,-1}^2x^2 - 1296(\alpha + 9\varkappa)b_{1,-1}x - 7776\varkappa)
\end{gathered}
\end{equation}
where $C_2$ is the integration constant to be determined by expanding the right-hand side of the
equation~\eqref{eq:A2-C2} in a Taylor series at $x=0$,
\begin{equation}\label{eq:A2-Taylor}
\begin{aligned}
\overset{\;\infty}{A}_2(x)&=-\frac{3\mi\varkappa}{b_{1,-1}}x-\frac{\mi}{2}(\alpha+12\varkappa)x^2\\
&-\left(\frac{C_2}{36}+\left(\frac{44\mi\varkappa}{3}+\frac{3\mi\alpha}{4}
+\frac{151\mi}{12}\right)b_{1,-1}\right)x^3+\mathcal{O}(x^4).
\end{aligned}
\end{equation}
Comparing the expansion~\eqref{eq:A2-Taylor} with the expansion~\eqref{eq:Ak-inftyDef} for $k=2$, one gets
\begin{gather}
\Theta a_{1,1}=-\frac{3\mi\varkappa}{b_{1,-1}}=-\frac{3\mi\varkappa}{a_{1,-1}\Theta},\quad
\Theta^2 a_{2,0}=-\frac{\mi}{2}(\alpha+12\varkappa),\label{eqs:a11-a20}\\
\Theta^3 a_{3,-1}=
-\frac{C_2}{36}-\left(\frac{44\mi\varkappa}{3}+\frac{3\mi\alpha}{4}+\frac{151\mi}{12}\right)b_{1,-1}.
\label{eq:C2a3-1}
\end{gather}
The first equation in \eqref{eqs:a11-a20} is equivalent to the last equation in \eqref{eq:coeffs:a1j}, the
second equation for $a_{2,0}$ is equivalent to the expression for this coefficient obtained in \cite{KitVar2023}
(cf. the last equation in (C.9) of \cite{KitVar2023}). The coefficient $a_{3,-1}$ was also computed in
\cite{KitVar2023} (see the last equation in (C.10) \cite{KitVar2023}), that allows one to find
\begin{equation}\label{eq:C2-formula}
C_2=\frac{3\mi}{2}(3\alpha^2+24\alpha\varkappa+30\varkappa^2-6\alpha-272\varkappa-305)b_{1,-1}.
\end{equation}
Substituting this value for $C_2$ into the equation~\eqref{eq:A2-C2} for $\overset{\;\infty}{A}_2(x)$ and
decomposing the resulting equation in partial fractions, we get
\begin{equation}\label{eq:A2-parfrac}
\begin{gathered}
\overset{\;\infty}{A}_2(x)=\frac{\mi}{b_{1,-1}^2}\left(-\frac{\varkappa+1}{12}b_{1,-1}^3x^3
-\frac{\alpha+12\varkappa+12}{2}b_{1,-1}^2x^2\right.\\
+\frac{3}{2}(\Xi_2-20\varkappa-53)b_{1,-1}x
\left.+36(\Xi_2+6\varkappa-23)+
\sum_{n=1}^3\frac{\eta_n}{(1-b_{1,-1}x/6)^n}\right),
\end{gathered}
\end{equation}
where
\begin{equation}\label{eq:Xi2}
\Xi_2=3(\alpha^2+8\alpha\varkappa+10\varkappa^2+4\alpha),
\end{equation}
\begin{equation}\label{eq:eta123}
\eta_1=-27(3\Xi_2+42\varkappa-47),\;
\eta_2=9(7\Xi_2+162\varkappa-59),\;
\eta_3=-18(\Xi_2+30\varkappa - 5).
\end{equation}
Expanding the right-hand side of the equation~\eqref{eq:A2-parfrac} into a Taylor series at $x=0$, one reproduces
the first three off-set coefficients: $a_{1,1}$, $a_{2,0}$, and $a_{3,-1}$ (cf. equation~\eqref{eq:A2-Taylor}
and a paragraph below), and, comparing the resulting expansion with the expansion~\eqref{eq:Ak-inftyDef}, arrives
at the general formula,
\begin{equation}\label{eq:an-(n-2)}
a_{n,-(n-2)}=\frac{\mi a_{1,-1}^{n-2}}{6^n\Theta^2}\left(\eta_1+(n+1)\eta_2+\frac{(n+1)(n+2)}{2}\eta_3\right),\qquad
n\geqslant 4.
\end{equation}
The formula~\eqref{eq:an-(n-2)} for $n=k$ is equivalent to the formula that was conjectured in \cite{KitVar2023}
(see Conjecture C.2 of \cite{KitVar2023}, for the lower choice of signs in equation (C.14)).

Now we will describe the continuation of the procedure for calculating the generating functions
$\overset{\;\infty}{A}_k(x)$. This procedure requires a priori knowledge of the coefficients $a_{k,-1}$,
$k\geqslant4$. The coefficients $a_{4,-1}=a_{6,-1}=0$ and $a_{5,-1}$ were calculated in \cite{KitVar2023}. At this stage we conjecture that
The functions $\overset{\;\infty}{A}_{2j-1}(x)\equiv0$ for all $j\in\mathbb{N}$. This statement can be proved by
mathematical induction provided that it is confirmed that $a_{2j,-1}=0$. Actually the induction hypothesis for
the function $\overset{\;\infty}{A}_1$ is established above, then, assuming that $\overset{\;\infty}{A}_{2i-1}(x)\equiv0$ for $i<j$,
one shows that $\overset{\;\infty}{A}_{2j-1}(x)$ is a solution of a second order linear homogeneous ODE, whose
rational solution reads,
\begin{equation}\label{eq:A2j-1-rational}
\overset{\;\infty}{A}_{2j-1}(x)=\frac{C_{2j-1}x^{2j}(xb_{1,-1}+6)}{(xb_{1,-1}-6)^3}=
-\frac{C_{2j-1}}{36}x^{2j}+\mathcal{O}(x^{2j+1}),
\end{equation}
where $C_{2j-1}$ is the integration constant. Comparing the expansion~\eqref{eq:A2j-1-rational} with the
expansion~\eqref{eq:Ak-inftyDef} for $k=2j-1$, we find that $-C_{2j-1}/36=a_{2j,-1}\Theta^{2j}$; this would complete
the inductive step provided it is established that $a_{2j,-1}=0$.\footnote{\label{foot:ak-1} One can prove that
$a_{2j,-1}=0$ and calculate $a_{2j-1,-1}$ for all $j\in\mathbb{N}$ with the help of the conjugate expansion constructed
in Subsection~\ref{subsec:trigInftyB}, see Remarks~\ref{rem:Bn-even} and \ref{rem:Cn-an|pm1}.}

Calculating the generating functions $\overset{\;\infty}{A}_{2j}(x)$ for $j\geqslant2$ is technically more
complicated, but practically resembles the case $j=1$.
It is helpful that \textsc{Maple} contains an implementation of an algorithm for constructing rational solutions
to linear ODEs, but these solutions quickly become quite complex as $j$ grows, so making the results of the computer
calculations comprehensible takes considerably more time than constructing the rational solution itself.
As an example, consider the calculation of $\overset{\;\infty}{A}_{4}(x)$, which takes less than 1 second on a
computer, but took the author several hours to find the representation of the results given below. It seems that
there may be more appropriate variables and relatively simple algebraic algorithms that could help to write explicit
formulae for the generating functions of $\overset{\;\infty}{A}_{2j}(x)$ in the general case.
\begin{equation}\label{eq:A4C4}
\begin{gathered}
\overset{\;\infty}{A}_{4}(x)=\frac{C_4x^5(xb_{1,-1}+6)}{(xb_{1,-1}-6)^3}
-\frac{x^2}{864b_{1,-1}^2(xb_{1,-1}-6)^4}\Big(2(\varkappa + 1)^2b_{1,-1}^8x^8\\
+3(\varkappa+1)(\Xi_2+24\alpha+124\varkappa+91)b_{1,-1}^7x^7\\
-72(\varkappa+1)(\Xi_2+24\alpha+134\varkappa+101)b_{1,-1}^6x^6\\
+54(\Xi_2^2+60\varkappa\Xi_2+52\Xi_2+336\alpha(2\varkappa+1)\!+5068\varkappa^2\!+5728\varkappa\!+1539)
b_{1,-1}^5x^5\\
+648(\Xi_2^2+60\varkappa\Xi_2-38\Xi_2-96\alpha(2\varkappa+1)-548\varkappa^2-2588\varkappa-423)b_{1,-1}^4x^4\\
+93312\varkappa(\Xi_2-24\alpha-74\varkappa-3)b_{1,-1}^2x^2\\
-139968\varkappa(\Xi_2-24\alpha-64\varkappa-3)+3359232\varkappa^2\Big),
\end{gathered}
\end{equation}
where $\Xi_2$ is defined in the equation~\eqref{eq:Xi2}. Expanding $\overset{\;\infty}{A}_{4}(x)$ given by
\eqref{eq:A4C4} in the Taylor series at $x=0$, one gets
\begin{equation}\label{eq: A4C4-Taylor}
\begin{gathered}
\overset{\;\infty}{A}_{4}(x)=-\frac{3\varkappa^2}{b_{1,-1}^2}x^2
+\frac{\varkappa}{8b_{1,-1}}(\Xi_2-24\alpha-80\varkappa-3)x^3\\
-\left(\varkappa\big(3\Xi_2-72\alpha-232\varkappa-9\big)\frac{b_{1,-1}}{144}+\frac{C_4}{36}\right)x^5
+\mathcal{O}(x^6).
\end{gathered}
\end{equation}
Comparing expansions~\eqref{eq: A4C4-Taylor} and \eqref{eq:Ak-inftyDef} for $k=4$, we find that
\begin{align}
a_{2,2}\Theta^2&=-\frac{3\varkappa^2}{b_{1,-1}^2}=-\frac{3\varkappa^2}{\Theta^2a_{1,-1}^2},\;\;
\Rightarrow\;\; a_{2,2}=a_{1,1}^2/3,\quad
a_{4,0}=0,\label{eq:a22-a40}\\
a_{3,1}\Theta^3&=\frac{\varkappa}{8b_{1,-1}}(\Xi_2-24\alpha-80\varkappa-3)\nonumber\\
&=\frac{\varkappa}{8a_{1,-1}\Theta}(3\alpha^2+24\varkappa\alpha+30\varkappa^2-12\alpha-80\varkappa-3),\label{eq:a31}\\
a_{5,-1}\Theta^5&=-\left(\varkappa\big(3\Xi_2-72\alpha-232\varkappa-9\big)\frac{b_{1,-1}}{144}+\frac{C_4}{36}\right)
\label{eq:a5m1}
\end{align}
The formulae for the coefficients $a_{2,2}$ and $a_{4,0}$ in \eqref{eq:a22-a40} coincide with the corresponding
formulae obtained in \cite{KitVar2023} (cf. with the equations (C.9) and (C.11) in \cite{KitVar2023},
respectively).
Taking into account that $\varkappa=\mi a_{1,1}a_{1,-1}\sqrt{3}(\varepsilon b)^{1/3}$
(see the last equation in~\eqref{eq:coeffs:a1j} and the first equation in~\eqref{eqs:varkappa-alpha},
and using the relation $\Theta^4=3^3(\varepsilon b)^{2/3}$ (see the last equation in \eqref{eq:xy-infty1}),
it is confirmed that the formulae for $a_{3,1}$ given by equation~\eqref{eq:a31} and the formula for $a_{3,1}$
presented in the paper \cite{KitVar2023} are equivalent (cf. with the equation (C.10) in \cite{KitVar2023}
with the choice of upper signs).
The coefficient $a_{5,-1}$ is also calculated in \cite{KitVar2023} (see equation (C.12) for the lower signs),
this allows us to find the constant of integration $C_4$,
\begin{equation}\label{eq:C4-explicit}
\begin{aligned}
C_4&=\frac{b_{1,-1}}{32}\Big(\Xi_2^2-16\alpha\Xi_2-72\varkappa\Xi_2-14\Xi_2\\
&+48\alpha(4\varkappa-1)-2000\varkappa^3-1892\varkappa^2- 168\varkappa+81\Big).
\end{aligned}
\end{equation}
Substituting the value of $C_4$, given by the equation~\eqref{eq:C4-explicit}, into the equation~\eqref{eq:A4C4} and
decomposing the resulting formula for $\overset{\;\infty}{A}_{4}(x)$ into partial fractions, one obtains
\begin{equation}\label{eq:A4-parfrac}
\begin{gathered}
\overset{\;\infty}{A}_{4}(x)=\!-\frac{1}{b_{1,-1}^4}\Bigg(\frac{(\varkappa + 1)^2}{432}b_{1,-1}^6x^6\!
+\!\frac{\varkappa+1}{288}\big(\Xi_2\!+\!24\alpha+\!140\varkappa+\!107\big)b_{1,-1}^5x^5\\
+\frac{1}{32}\big(\Xi_2^2+(16\alpha+168\varkappa+94)\Xi_2+144\alpha(4\varkappa+1)
+2000\varkappa^3+8732\varkappa^2\\
+5824\varkappa+493\big)b_{1,-1}^3x^3
+\frac{3}{2}\big(\Xi_2^2+(8\alpha+116\varkappa+30)\Xi_2+24\alpha(10\varkappa+3)\\
+1000\varkappa^3\!+\!4370\varkappa^2\!+\!2108\varkappa\!+\!247\big)b_{1,-1}x^2\!
+\!\frac{81}{8}\big(3\Xi_2^2+\!(16\alpha+292\varkappa+42)\Xi_2\\
+48\alpha(10\varkappa+3)+2000\varkappa^3+9444\varkappa^2+3484\varkappa+511\big)b_{1,-1}x
+432\big(\Xi_2^2\\
+(4\alpha+88\varkappa+7)\Xi_2+12\alpha(10\varkappa+3)+
+500\varkappa^3+2562\varkappa^2+742\varkappa+133\big)\\
-\sum_{n=1}^4\frac{\nu_n}{(1-b_{1,-1}x/6)^n}\Bigg),
\end{gathered}
\end{equation}
where
\begin{align*}
\nu_1&=\frac{135}{4}(25\Xi_2^2+(80\alpha+2060\varkappa+78)\Xi_2+240\alpha(10\varkappa+3)\nonumber\\
&+10000\varkappa^3+55452\varkappa^2+12692\varkappa+2773),\\
\nu_2&=-\frac{27}{4}(91\Xi_2^2+(176\alpha+6692\varkappa-246)\Xi_2+528\alpha(10\varkappa+3)\nonumber\\
&+22000\varkappa^3+153012\varkappa^2+14012\varkappa+6943),\\
\nu_3&=27(9\Xi_2^2+(8\alpha+596\varkappa-62)\Xi_2+24\alpha(10\varkappa+3)\\
&+5(200\varkappa^3+2256\varkappa^2-188\varkappa+87)),\\
\nu_4&=-\frac{81}{2}(\Xi_2+30\varkappa-5)^2.
\end{align*}
Expanding now $\overset{\;\infty}{A}_{4}(x)$ given by the equation~\eqref{eq:A4-parfrac} into a Taylor series
at $x=0$, we reproduce the first four off-set coefficients presented by
equations~\eqref{eq:a22-a40}--\eqref{eq:a5m1}; note that the coefficients $\nu_k$ satisfy the relation
\begin{equation}\label{eq:nu-relation}
\nu_{1}+5\nu_{2}+15\nu_{3}+35\nu_{4}=0.
\end{equation}
from which it follows that the coefficient $a_{4,0}=0$. Comparing the coefficients of the Taylor expansion with
the corresponding coefficients of the expansion defining $\overset{\;\infty}{A}_{k}(x)$
(see equation~\eqref{eq:Ak-inftyDef} for $k=4$, we get at the power $x^6$ the last off-set coefficient,
\begin{equation}\label{eq:a6m2}
\begin{aligned}
a_{6,-2}=&-\frac{a_{1,-1}^2}{864\Theta^4}\Big(\Xi_2^2-(8\alpha-4\varkappa+26)\Xi_2\\
&-24\alpha(10\varkappa+3)-1000\varkappa^3-1990\varkappa^2-1408\varkappa-171\Big).
\end{aligned}
\end{equation}
The formula~\eqref{eq:a6m2} coincides with the formula for $a_{6,-2}$ obtained in \cite{KitVar2023} (cf. the
last equation in the list (C.13) \cite{KitVar2023} for the lower choice of signs). For the remaining coefficients
of the expansion~\eqref{eq:Ak-inftyDef} for $k=4$, we get the general formula
\begin{equation}\label{eq:an4-n}
a_{n,4-n}=-\frac{a_{1,-1}^{n-4}}{6^n\Theta^4}\sum_{l=1}^4\binom{n+l-1}{n}\nu_l,\qquad
n\geqslant7.
\end{equation}
\begin{remark}\label{rem:symmetry-infty}
The  expansion~\eqref{eq:app:u-asympt-regular-complete} is uniquely determined
by its first non-trivial coefficients $a_{1,1}$ and $a_{1,-1}$ and, of course, by the coefficients
of the equation~\eqref{eq:dp3u}. In particular, this fact implies that the following symmetry property for the
coefficients $a_{k,j}$: any relation between the coefficients $a_{k,j}$ will remain valid if we simultaneously
replace all the second subscripts $j$ of the coefficients involved in this relation by $-j$, as well as the symbols
$\mi$, $\alpha$ and $\varkappa$, if they are present in the relation, by $-\mi$, $-\alpha$ and $-\varkappa$,
respectively. For example, applying this symmetry to the equation~\eqref{eq:ak-k}, one proves that
\begin{equation}\label{eq:ann}
a_{n,n}=\frac{na_{1,1}^n}{6^{n-1}}.
\end{equation}
According to the symmetry principle formulated above the equations~\eqref{eqs:a11-a20} are invariant under application
of the symmetry transformation. As a less elementary example of the application of this symmetry principle, we can
obtain the explicit formula for $a_{3,-1}$ with the help of the equation~\eqref{eq:a31},
\begin{equation}\label{eq:a3m1}
\begin{aligned}
a_{3,-1}&=-\frac{\varkappa}{8a_{1,1}\Theta^4}(3\alpha^2+24\varkappa\alpha+30\varkappa^2+12\alpha+80\varkappa-3)\\
&=-\frac{\mi \sqrt{3}a_{1,-1}}{6^3(\varepsilon b)^{1/3}}
(3\alpha^2+24\alpha\varkappa+30\varkappa^2+12\alpha+80\varkappa-3),
\end{aligned}
\end{equation}
which coincides with the corresponding formula (C.10) \cite{KitVar2023}.

Note that instead of variables $x(\tau)$ and $y(\tau)$ (cf. equation~\eqref{eq:xy-infty1}), one can define
symmetric variables $\tilde{x}(\tau)$ and $\tilde{y}(\tau)$,
\begin{equation}\label{eq:tide-x:tilde-y}
\begin{aligned}
\tilde{x}(\tau)&:=\theta(\tau)^{-1/2}\tilde{w}^{-1}=\Theta^{-1}\tau^{-1/3}\tau^{2\varkappa/3}\me^{\mi\theta(\tau)},\\
\tilde{y}(\tau)&:=\tilde{w}=w^{-1}=\tau^{-2\varkappa/3}\me^{-\mi\theta(\tau)},\\
\end{aligned}
\end{equation}
and reproduce a similar construction
for another super-generating function for infinite sequences of coefficients $a_{k,j}$, but with $j\to+\infty$.
\hfill$\blacksquare$\end{remark}
\subsection{Conjugate Expansion}\label{subsec:trigInftyB}
Let us consider the construction of the expansion conjugate to the expansion~\eqref{eq:A-inftyDef}.
Note that the equation~\eqref{eq:A-inftyPDE} is invariant under successive application of the following two
transformations, (i) $x\to y$, $\varkappa\to\tilde{\varkappa}-1/2$ and (ii)
$(\mi,\alpha,\tilde{\varkappa})\to(-\mi,-\alpha,-\tilde{\varkappa})$ (see Remark~\ref{rem:symmetry-infty}).
The shift of $\varkappa$ by $1/2$ should not be confusing, since the definition of the variables $x$ and $y$ is
not completely symmetric (the factor $\tau^{-1/3}$ in the definition~\eqref{eq:xy-infty1}). Thus,  if we substitute
into the equation~\eqref{eq:A-inftyPDE} for $\overset{\;\infty}{A}$ the expansion~\eqref{eq:AxyB} and
choose the general rational  solution for $\overset{\;\infty}{B}_0(y)=\overset{\;\infty}{A}_0(x)|_{x=y}$, then
the remaining generating functions $\overset{\;\infty}{B}_k(y)$, $k\geqslant1$ would coincide with the corresponding
functions $\overset{\;\infty}{A}_k(x)|_{x=y, \varkappa=\tilde{\varkappa}}$. The supergenerating function constructed
by this method for the functions $\overset{\;\infty}{B}_k(y)$ is equivalent to the supergenerating function
$\overset{\;\infty}{A}(x,y)$: the parameter $\tilde{\varkappa}$ has the same meaning as the original parameter
$\varkappa$, it is just a difference in the notation of this parameter as well as the substitution
$x\leftrightarrow y$, however, these super-generating functions are different, because the new supergenerating
function is an even function of $x$, while the function $\overset{\;\infty}{A}(x,y)$ is an even function of $y$.
Therefore, the construction of the above proposed ``$B$-expansion'' is not conjugate to the ``$A$-expansion''.

The conjugate expansion for $\overset{\;\infty}{A}(x,y)$ is constructed by taking the cut-off rational solution,
\begin{equation}\label{eq:B0}
\overset{\;\infty}{B}_0(y)=1.
\end{equation}
For $\overset{\;\infty}{B}_1(y)$ one gets the linear ODE,
$y\frac{\md^2}{\md y^2}\overset{\;\infty}{B}_1(y)-\frac{\md}{\md y}\overset{\;\infty}{B}_1(y)=0$, so that
\begin{equation}\label{eq:B1}
\overset{\;\infty}{B}_1(y)=b_{1,1}y^2+b_{1,-1},
\end{equation}
where $b_{1,\pm1}\in\mathbb{C}$ coincide with the corresponding parameters introduced in the previous subsection.
Equation for $\overset{\;\infty}{B}_2(y)$ reads
\begin{equation}\label{eq:diffeqB2}
2y^2\frac{\md^2}{\md y^2}\overset{\;\infty}{B}_2(y)-6y\frac{\md}{\md y}\overset{\;\infty}{B}_2(y)+
6\overset{\;\infty}{B}_2(y)=2b_{1,1}^2y^4+(\mi\alpha-4b_{1,1}b_{1,-1})y^2+2b_{1,-1}^2,
\end{equation}
so that the general solution is a polynomial
\begin{equation}\label{eq:B2}
\overset{\;\infty}{B}_2(y)=C_4y^3+C_3y+\frac{b_{1,1}^2}{3}y^4-\frac12(\mi\alpha-4b_{1,1}b_{1,-1})y^2+
\frac{b_{1,-1}^2}{3},
\end{equation}
where $C_3,C_4\in\mathbb{C}$ are integration constants that need to be determined. Continuing this procedure
one obtains for the function $\overset{\;\infty}{B}_3(y)$ the linear second order ODE,
\begin{gather*}
3y^2\frac{\md^2}{\md y^2}\overset{\;\infty}{B}_3(y)-15y\frac{\md}{\md y}\overset{\;\infty}{B}_3(y)
+24\overset{\;\infty}{B}_3(y)=\nonumber\\
2b_{1,1}^3y^6+2(3\mi\varkappa+b_{1,1}b_{1,-1})(b_{1,1}y^2+b_{1,-1})y^2\\
+2b_{1,-1}^3+6(C_4y^2-C_3)(b_{1,1}y^2-b_{1,-1})y,
\end{gather*}
that has the following general solution,
\begin{equation}\label{eq:B3ln}
\begin{gathered}
\overset{\;\infty}{B}_3(y)=C_6y^4+C_5y^2+\frac{1}{3}(3\mi\varkappa+b_{1,1}b_{1,-1})(b_{1,1}y^2-b_{1,-1})y^2\ln(y)\\
+\frac{b_{1,1}^3}{12}y^6-\frac{1}{6}(3\mi\varkappa+b_{1,1}b_{1,-1})(b_{1,1}y^2+b_{1,-1})y^2+\frac{b_{1,-1}^3}{12}\\
+\frac{2C_4b_{1,1}}{3}y^5+2(C_4b_{1,-1}+C_3b_{1,1})y^3 +\frac{2C_3b_{1,-1}}{3}y,
\end{gathered}
\end{equation}
where $C_5,C_6\in\mathbb{C}$ are integration constants. The rationality condition leads to the relation
\begin{equation}\label{eq:varkappab11b1m1}
3\mi\varkappa+b_{1,1}b_{1,-1}=0,
\end{equation}
so that
\begin{equation}\label{eq:B3-no-ln}
\begin{gathered}
\overset{\;\infty}{B}_3(y)=\frac{b_{1,1}^3}{12}y^6+\frac{2}{3}C_4b_{1,1}y^5+C_6y^4\\
+2(C_4b_{1,-1}+C_3b_{1,1})y^3+C_5y^2+\frac{2}{3}C_3b_{1,-1}y+\frac{b_{1,-1}^3}{12}.
\end{gathered}
\end{equation}
At this stage, the function $\overset{\;\infty}{B}_2(y)$ depends on two complex parameters, $C_3$ and $C_4$, and
the function $\overset{\;\infty}{B}_3(y)$ depends on four complex parameters, $C_3$, $C_4$, $C_5$, and $C_6$.
Let us describe how the algorithm for determining these parameters works. In general, the linear ODE for the function
$\overset{\;\infty}{B}_n(y)$ for $n\geqslant4$ is
\begin{equation}\label{eq:eqdiff-Bn-general}
\begin{gathered}
y^2\frac{\md^2}{\md y^2}\overset{\;\infty}{B}_n(y)-(2n-1)y\frac{\md}{\md y}\overset{\;\infty}{B}_n(y)+
(n^2-1)\overset{\;\infty}{B}_n(y)\\
=P_{2n}(y; b_{1,-1},b_{1,1},C_{2n-5},C_{2n-4},C_{2n-3},C_{2n-2}),
\end{gathered}
\end{equation}
where $P_{2n}(y;b_{1,-1},b_{1,1},C_{2n-5},C_{2n-4},C_{2n-3},C_{2n-2})$ is a polynomial in $y$ of degree $2n$,
a polynomial of degree $n$ in two variables $b_{1,-1}$ and $b_{1,1}$,
and a linear function of parameters $C_k$, $k=2n-5,2n-4,2n-3,2n-2$,\footnote{\label{foot:quadraticC34B4}
In fact, this statement is not true for the polynomial $P_8(y;\ldots)$ and, correspondingly, for the polynomial
$Q_8(y;\ldots)$ (see equation~\eqref{eq:Bn-general} for $n=4$ below), the coefficients of the terms $y^6$ and $y^2$
of these polynomials contain quadratic dependencies $C_3^2$ and $C_4^2$, respectively. This circumstance does not
create any problems when constructing a supergenerating function.}
which appeared during the integration of the corresponding linear ODEs for the functions
$\overset{\;\infty}{B}_{n-2}$ and $\overset{\;\infty}{B}_{n-1}$
and have not yet been determined.
The general solution of the homogeneous part of the equation~\eqref{eq:eqdiff-Bn-general},
$\overset{\;\infty}{B}_n^{hom}(y)$, reads,
\begin{equation}\label{eq:BnHOM}
\overset{\;\infty}{B}_n^{hom}=C_{2n}y^{n+1}+C_{2n-1}y^{n-1}.
\end{equation}
The general solution of the equation~\eqref{eq:eqdiff-Bn-general} can be written as
\begin{equation}\label{eq:Bn-general}
\begin{aligned}
\overset{\;\infty}{B}_n^{gen}(y)=\overset{\;\infty}{B}_n^{hom}(y)
&+Q_{2n}(y;b_{1,-1},b_{1,1},C_{2n-5},C_{2n-4},C_{2n-3},C_{2n-2})\\
&+L_1(b_{1,-1},b_{1,1},C_{2n-5},C_{2n-4})y^{n+1}\ln(y)\\
&+L_2(b_{1,-1},b_{1,1},C_{2n-5},C_{2n-4})y^{n-1}\ln(y),
\end{aligned}
\end{equation}
where $Q_{2n}(y;b_{1,-1},b_{1,1},C_{2n-5},C_{2n-4},C_{2n-3},C_{2n-2})$, is a polynomial function of the same type as
$P_{2n}(y;b_{1,-1},b_{1,1},C_{2n-5},C_{2n-4},C_{2n-3},C_{2n-2})$,\\ $L_k(b_{1,-1},b_{1,1},C_{2n-5},C_{2n-4})$, $k=1,2$,
are linear functions of $C_{2n-5}$ and $C_{2n-4}$. For brevity, we denote
$L_k:=L_k(b_{1,-1},b_{1,1},C_{2n-5},C_{2n-4}$.

For $n=2k+2$, $k\in\mathbb{N}$, one proves that
\begin{equation}\label{eqs:L1L2-n-even}
\begin{aligned}
L_1&=\frac{1}{6}\left(2b_{1,1}^2C_{4k-1}+\big(2b_{1,1}b_{1,-1}-3\mi(2k-1)\big)C_{4k}\right),\\
L_2&=-\frac{1}{6}\left(2b_{1,-1}^2C_{4k}+\big(2b_{1,1}b_{1,-1}+3\mi(2k-1)\big)C_{4k-1}\right).
\end{aligned}
\end{equation}
For $n=2k+3$, $k\in\mathbb{N}$, one finds that
\begin{equation}\label{eqs:L1L2-n-odd}
\begin{aligned}
L_1&=\frac{1}{3}\left(b_{1,1}^2C_{4k+1}+\big(b_{1,1}b_{1,-1}-3\mi k\big)C_{4k+2}\right)
+b_{1,1}F_{2k}^1(\varkappa,\alpha),\\
L_2&=-\frac{1}{3}\left(b_{1,-1}^2C_{4k+2}+\big(b_{1,1}b_{1,-1}+3\mi k\big)C_{4k+1}\right)
+b_{1,-1}F_{2k}^2(\varkappa,\alpha),
\end{aligned}
\end{equation}
where $\varkappa$ is defined in equation~\eqref{eq:varkappab11b1m1}, and $F_{2k}^l(\varkappa,\alpha)$, $l=1,2$,
are polynomials in two variables $\varkappa$ and $\alpha$ of degree $2k$.

Now we can describe the algorithm for constructing the polynomials $\overset{\;\infty}{B}_n(y)$. The function
$\overset{\;\infty}{B}_n^{gen}(y)$ depends
on three pairs of complex parameters: $C_{2n},C_{2n-1}$, $C_{2n-2},C_{2n-3}$, and $C_{2n-4},C_{2n-5}$,
which are to be determined. By equating the coefficients of the logarithmic terms to zero
(see equation~\eqref{eq:Bn-general}), we find that the pair of parameters with the smallest subscripts,
$C_{2n-4},C_{2n-5}$, are solutions to the system of two linear equations, $L_1=0$ and $L_2=0$, with non-zero
determinant, so that they can be determined uniquely. Thus, we arrive at a polynomial $\overset{\;\infty}{B}_n(y)$
that depends on two pairs of undetermined parameters: $C_{2n},C_{2n-1}$ and $C_{2n-2},C_{2n-3}$. This fact, in turn,
means that polynomial $\overset{\;\infty}{B}_{n-1}(y)$ depends on one pair of parameters $C_{2n-2},C_{2n-3}$ and
the polynomial $\overset{\;\infty}{B}_{n-2}(y)$ is completely determined. Therefore, by a mathematical induction
argument, one proves the existence and uniqueness of the sequence of polynomials
$\overset{\;\infty}{B}_n(y)$, $n\in\mathbb{N}$.

Two important remarks regarding the construction of the polynomials $\overset{\;\infty}{B}_n(y)$ should be made.
\begin{remark}\label{rem:Bn-even}
The equations~\eqref{eqs:L1L2-n-even} are homogeneous, so that $C_{4k-1}=C_{4k}=0$ for $k\in\mathbb{N}$, since,
non-homogeneous terms in these equations can appear only as linear combinations of the coefficients of products of
monomials of the $\overset{\;\infty}{B}_m(y)$ polynomials with $m<n$ that contain an odd number of monomials of
odd degree, but such monomials are absent in the $\overset{\;\infty}{B}_m(y)$ polynomials, which can be justified
using the mathematical induction argument. This observation confirms that
$\overset{\;\infty}{A}_{2j-1}(x)=0$ (cf. the paragraph above equation~\eqref{eq:A2j-1-rational}).
\hfill$\blacksquare$\end{remark}
\begin{remark}\label{rem:Cn-an|pm1}
\begin{equation}\label{eq:Cn-an|pm1}
a_{n,1}=C_{2n}\Theta^{-n},
\quad
a_{n,-1}=C_{2n-1}\Theta^{-n},
\qquad
n\geqslant2.
\end{equation}
Thus as follows from Remark~\ref{rem:Bn-even}, $a_{2k,1}=a_{2k,-1}=0$ for $k\in\mathbb{N}$. Putting $n=2k-1$,
we find that the coefficients $a_{2k+1,\pm1}$, which the key coefficients for constructing the generating functions
$\overset{\;\infty}{A}_{2k}(x)$, are computed via the coefficients $C_{4k+2}$ and $C_{4k+1}$.
\hfill$\blacksquare$\end{remark}
\subsection{Application to Asymptotics of General Solutions: $\varkappa\neq0$}\label{subsec:genInftyAsympt}
The goal of this subsection is to show that the super-generating function
constructed in the previous subsection, under a proper restrictions on the parameters, can serve as an asymptotic
expansion for the solution $u(\tau)$. Of course, we started Subsection~\ref{subsec:Ainfty:trig} with the asymptotic
expansion for the general solutions of the equation~\eqref{eq:dp3u}, however, here we obtain an analytic continuation
of the latter expansion to the new values of its parameters.

In this section we consider the expansion constructed in Subsection~\ref{subsec:Ainfty:trig} but for the
symmetric variables $\tilde{x}(\tau)$ and $\tilde{y}(\tau)$ defined in Remark~\ref{rem:symmetry-infty}. All
formulae constructed in the previous subsection for variables $x(\tau)$ and $y(\tau)$ remain valid for
variables $\tilde{x}(\tau)$ and $\tilde{y}(\tau)$ after the symmetry transformation described in
Remark~\ref{rem:symmetry-infty}.
We turn to the symmetric variables because we would like to get exact correspondence
with the asymptotic formulae obtained in \cite{KitVar2023}.\footnote{\label{foot:asympt-infty} The formulae in terms
of variables $x(\tau)$ and $y(\tau)$ also can be used to define correct asymptotics, which is valid  in a different
domain of the parameters (see Subsection~\ref{subsec:truncatedsolutions}.)}

Consider the symmetric (in the sense of Remark~\ref{rem:symmetry-infty}) to $\overset{\;\infty}{A}(x,y)$
super-generating function
\begin{equation}\label{eq:tildeAinfty-trig}
\stackrel{\;\infty}{\tilde{A}}(\tilde{x},\tilde{y})=\sum_{k=0}^{\infty}\tilde{y}^k
\stackrel{\;\infty}{\tilde{A}}_k(\tilde{x}),
\end{equation}
and substitute in it $\tilde{x}=\tilde{x}(\tau)$ and $\tilde{y}=\tilde{y}(\tau)$ defined by
the equations~\eqref{eq:tide-x:tilde-y} as a result we have an expansion that depend on $\tau$ which we are
going to treat as an asymptotic expansion for the function $u(\tau)$ in the same domain $\Omega$
(see the equation\eqref{eq:Omega-infty-def}) as the original asymptotic
expansion~\eqref{eq:app:u-asympt-regular-complete}.

The leading term of the, so far, proposed asymptotic expansion is
\begin{gather}
\overset{\;\infty}{\tilde{A}}_0(\tilde{x}(\tau))=1+
\frac{b_{1,1}\tilde{x}(\tau)}{\left(1-b_{1,1}\tilde{x}(\tau)/6\right)^2}=
1-\frac{3}{2\sin^2(\frac12\vartheta(\tau))},\label{eq:tildeAinfty-A0}\\
\vartheta(\tau)=\theta(\tau)-\mi(\varkappa-1/2)\ln\theta(\tau)-\mi\ln\big(a_{1,1}/6\big)
+\mi(\varkappa-1/2)\ln(\Theta^2),\label{eq:tildeAinfty-vartheta}
\end{gather}
where we used the symmetry defined in Remark~\ref{rem:symmetry-infty} to the equation~\eqref{eq:A0-inftyRat}.
Note that $b_{1,1}=\Theta a_{1,1}$, as follows from applying symmetry to the equation~\eqref{eq:ak-k} for $n=1$.

Formulae~\eqref{eq:tildeAinfty-A0} and \eqref{eq:tildeAinfty-vartheta} coincide, up to parameterization
$a_{1,1}$ by monodromy data, with the corresponding formulae (C.26) and (C.28) from \cite{KitVar2023}, which define
the leading term of the asymptotics $u(\tau)$ in the region $0<\mathrm{Re}\,(\varkappa)<1$.\footnote{\label{foot:a11}
The parametrization of $a_{1,1}$ via monodromy data is given by equation~(C.8) in \cite{KitVar2023}.}
At this stage, we have no restrictions on the parameter $\varkappa$; we will find them below, based on
the requirement that the expansion~\eqref{eq:tildeAinfty-trig} defines an asymptotic expansion for $u(\tau)$.

Consider first the case when $\tilde{x}(\tau)\to0$ (or $\mathcal{O}(1)$) as $\tau\to\infty$, then the first
equation~\eqref{eq:tide-x:tilde-y} implies that $\mathrm{Re}\,\varkappa\leqslant1/2$.
If $\mathrm{Re}\,\varkappa<1/2$, then $\overset{\;\infty}{\tilde{A}}_{2k}(\tilde{x})=\mathcal{O}(\tilde{x}^{k})$
as $\tilde{x}\to0$; this fact for $k=0$, $1$, and $2$ follows from the equations~\eqref{eq:A0-inftyRat},
\eqref{eq:A2-Taylor}, and \eqref{eq: A4C4-Taylor}, respectively.
If $\mathrm{Re}\,\varkappa=1/2$ ($\tilde{x}=\mathcal{O}(1)$), then the equations~\eqref{eq:A0-inftyRat},
\eqref{eq:A2-parfrac}, and \eqref{eq:A4-parfrac} imply that
$\overset{\;\infty}{\tilde{A}}_{2k}(\tilde{x}(\tau))$, $k=0,1,2$  are bounded functions
of $\tau$ in $\Omega$.\footnote{\label{foot:Ainfty2k-1=0} Recall that
$\overset{\;\infty}{\tilde{A}}_{2k-1}(\tilde{x}(\tau))=0$, $k=\in\mathbb{N}$, see Remark~\ref{rem:Bn-even}.}
\hspace{-8pt}${}^,$\footnote{\label{foot:cheese-likeOmega}If $\varkappa=1/2$, then for a proper choice of $a_{1,1}$ there might be
an infinite sequence $\tau_n$ of zeros, $\vartheta(\tau_n)=2\pi n$ of the function $\sin$ in the denominator
of \eqref{eq:tildeAinfty-A0}. For solutions corresponding to this choice of $a_{1,1}$ one has to redefine the
domain $\Omega$ (see the equation~\eqref{eq:Omega-infty-def}) as
$\Omega^*=\Omega\setminus\underset{n}{\cup}C_n(\delta)$, $C_n(\delta)=|\vartheta(\tau)-\vartheta(\tau_n)|<\delta<c$,
which is called sometimes a cheese-like domain.}
To confirm the above mentioned behaviour of $\overset{\;\infty}{\tilde{A}}_{2k}(\tilde{x}(\tau))$
for $k\in\mathbb{N}$, one has to apply a mathematical induction argument to a recurrence relation for the sequence
$\overset{\;\infty}{\tilde{A}}_{2k}(\tilde{x})$, which can be obtained by substituting the
expansion~\eqref{eq:A-inftyDef} into the partial differential equation~\eqref{eq:A-inftyPDE}.
Now we can obtain a condition that
the expansion~\eqref{eq:tildeAinfty-trig} after substituting $\tilde{x}(\tau)$ and $\tilde{y}(\tau)$ for $\tilde{x}$
and $\tilde{y}$, respectively, define an asymptotic series,
\begin{equation}\label{eq:Ainfty-trig-asymptestimate-smallx}
\begin{aligned}
\frac{(\tilde{y}(\tau))^{2k+2}\overset{\;\infty}{\tilde{A}}_{2k+2}(\tilde{x}(\tau))}
{(\tilde{y}(\tau))^{2k}\overset{\;\infty}{\tilde{A}}_{2k}(\tilde{x}(\tau))}
&=\mathcal{O}\left((\tilde{y}(\tau))^2\tilde{x}(\tau)\right)\\
&=\mathcal{O}\left(\tau^{-1/3-2\varkappa/3}\right)
\quad\Rightarrow\quad
\mathrm{Re}\,\varkappa>-\frac12.
\end{aligned}
\end{equation}
Thus, we see
that the expansion~\eqref{eq:tildeAinfty-trig} for $\tilde{x}=\tilde{x}(\tau)$ and $\tilde{y}=\tilde{y}(\tau)$
for $\mathrm{Re}\,\varkappa\in(-1/2,1/2]$ defines an asymptotic series. We impose one more requirement, namely,
the term with $\sin^2\big(\tfrac12\vartheta(\tau)\big)$ in the denominator of the second
equation~\eqref{eq:tildeAinfty-A0} is the leading term of the asymptotic series, i.e., it should be ``visible''
against the second term of the series, $(\tilde{y}(\tau))^{2}\overset{\;\infty}{\tilde{A}}_{2}(\tilde{x}(\tau))$,
\begin{equation}\label{eq:Ainf-trig-leading-cond-small}
(\tilde{y}(\tau))^{2}\overset{\;\infty}{\tilde{A}}_{2}(\tilde{x}(\tau))=o\big(\tilde{x}(\tau)\big)
\quad\Rightarrow\quad
\mathcal{O}\big((\tilde{y}(\tau))^{2}\big)=o(1)\quad\Rightarrow\quad
\mathrm{Re}\,\varkappa>0.
\end{equation}

Consider the case when $\tilde{x}(\tau)$
is growing. In terms of the parameter $\varkappa$ it means that $\mathrm{Re}\,\varkappa>1/2$. In this case,
$\overset{\;\infty}{\tilde{A}}_{2k}(\tilde{x}(\tau))=\mathcal{O}((x(\tau))^{3k})$. For $k=0$, $1$, and $2$, it follows
from the explicit formulae~\eqref{eq:A0-inftyRat}, \eqref{eq:A2-parfrac}, and \eqref{eq:A4-parfrac}, respectively.
For general $k\in\mathbb{N}$ can be confirmed by the mathematical induction argument, we left this proof for the
interested reader. Thus, the requirement that the expansion~\eqref{eq:tildeAinfty-trig} after substituting
$\tilde{x}=\tilde{x}(\tau)$ and $\tilde{y}=\tilde{y}(\tau)$ defines an asymptotic series is,
\begin{equation}\label{eq:Ainfty-trig-asymptestimate}
\begin{aligned}
\frac{(\tilde{y}(\tau))^{2k+2}\overset{\;\infty}{\tilde{A}}_{2k+2}(\tilde{x}(\tau))}
{(\tilde{y}(\tau))^{2k}\overset{\;\infty}{\tilde{A}}_{2k}(\tilde{x}(\tau))}
&=\mathcal{O}\left((\tilde{y}(\tau))^2(\tilde{x}(\tau))^3\right)\\
&=\mathcal{O}\left(\tau^{-1/3+2(\varkappa-1)/3}\right)
\quad\Rightarrow\quad
\mathrm{Re}\,\varkappa<\frac32.
\end{aligned}
\end{equation}
If we want that the term with $\sin^2\big(\tfrac12\vartheta(\tau)\big)$ in the denominator of the second
equation~\eqref{eq:tildeAinfty-A0} to determine the leading term of the asymptotic series, then this condition
should be strengthen; the visibility of this term against the second term of the asymptotic series means that
\begin{gather}
(\tilde{y}(\tau))^{2}\overset{\;\infty}{\tilde{A}}_{2}(\tilde{x}(\tau))=o\big(1/\tilde{x}(\tau)\big)
\quad\Rightarrow\quad
(\tilde{y}(\tau))^{2}(\tilde{x}(\tau))^{4}=o(1)
\quad\Rightarrow\quad\label{eq:Ainf-trig-leading-cond}\\
-\frac23-\frac23+\frac43\mathrm{Re}\,\varkappa<0,
\quad\Rightarrow\quad
\mathrm{Re}\,\varkappa<1.
\end{gather}
We can summarize the above study as follows, the formal expansion~\eqref{eq:tildeAinfty-trig} after substituting
$\tilde{x}=\tilde{x}(\tau)$ and $\tilde{y}=\tilde{y}(\tau)$ defines an asymptotic series for
$\mathrm{Re}\,\varkappa\in(-1/2,3/2)$. For $\mathrm{Re}\,\varkappa\in(-1/2,1/2)$ this series (after the re-expansion
of the functions $\overset{\;\infty}{\tilde{A}}_{2k}(\tilde{x}(\tau))$ coincide with the original asymptotic
series~\eqref{eq:app:u-asympt-regular-complete}, in that sense both series are equivalent for these values of
$\varkappa$. On the practical viewpoint the asymptotic series~\eqref{eq:app:u-asympt-regular-complete} gives better
approximation for the solutions, then the asymptotic series defined by the expansion~\eqref{eq:tildeAinfty-trig}
for $\mathrm{Re}\,\varkappa\in(-1/6,1/6)$, while the latter series ``works'' better  for
$\mathrm{Re}\,\varkappa\in(1/6,1/2)$. The length of these segments are conventional, of course, since it depends
also on the values of the other parameters. Thus, asymptotic series defined by
the expansion~\eqref{eq:tildeAinfty-trig} can be viewed as an analytic continuation of the asymptotic
series~\eqref{eq:app:u-asympt-regular-complete} with respect to the parameter $\varkappa$ to the strip
$\mathrm{Re}\,\varkappa\in[1/2,3/2)$. Furthermore, if $\mathrm{Re}\,\varkappa\in(0,1)$, then
$\overset{\;\infty}{\tilde{A}}_{0}(\tilde{x}(\tau))$ given by the equation~\eqref{eq:tildeAinfty-A0}
is the sum of two largest terms of the asymptotic series. In \cite{KitVar2010} it was proved that for any
given $\varkappa$ with $\mathrm{Re}\,\varkappa\in(0,1)$ and $a_{1,1}\in\mathbb{C}$, there exists
a unique solution $u(\tau)$ such leading term of its asymptotics coincide with
$\overset{\;\infty}{\tilde{A}}_{0}(\tilde{x}(\tau))$. Asymptotic series
$\overset{\;\infty}{\tilde{A}}(\tilde{x}(\tau),\tilde{y}(\tau))$, represent full asymptotic series of this
solution because all its terms are uniquely determined in terms of
$\overset{\;\infty}{\tilde{A}}_{0}(\tilde{x}(\tau))$  by substituting into the equation~\eqref{eq:dp3u}.

The function $\overset{\;\infty}{A}(x,y)$ after substituting $x=x(\tau)$ and $y=y(\tau)$, where
the functions $x(\tau)$ and $y(\tau)$ are defined by the equations~\eqref{eq:xy-infty1} define analytic continuation
with respect to the parameter $\varkappa$ to the strip $-3/2<\mathrm{Re}\,\varkappa\leqslant-1/2$ of the asymptotic
series~\eqref{eq:app:u-asympt-regular-complete}. This asymptotic series is equivalent to the asymptotic series
defined by $\overset{\;\infty}{\tilde{A}}(\tilde{x},\tilde{y})$; both asymptotic series define the same solution
$u(\tau)$ of the equation~\eqref{eq:dp3u}.

One can define (two equivalent) asymptotic series based on the conjugate expansion constructed in
Subsection~\ref{subsec:trigInftyB}; both asymptotics represent the same asymptotic
series~\eqref{eq:app:u-asympt-regular-complete}.

The supergenerating functions $\overset{\;\infty}{\tilde{A}}$ and $\overset{\;\infty}{\tilde{B}}$ allow, instead of
three different asymptotic expansions of the solution $u(\tau)$ for different ranges of the parameter
$\mathrm{Re}\,\varkappa$, to define a uniform asymptotic expansion valid for $\mathrm{Re}\,\varkappa\in(-3/2,3/2)$,
\begin{equation}\label{eq:u-trigAsympt-symmetric}
\begin{gathered}
u(\tau)=\frac{\varepsilon(\varepsilon b)^{2/3}}{2}\tau^{1/3}\left(\sum_{k=0}^N
(y(\tau))^k\overset{\;\infty}A_{k}(x(\tau))+
(\tilde{y}(\tau))^k\overset{\;\infty}{\tilde{A}}_{k}(\tilde{x}(\tau))\right.\\
-(x(\tau))^k\overset{\;\infty}B_{k}(y(\tau))
+\left.\mathcal{O}\Big((x(\tau))^{N+1}\overset{\;\infty}B_{N+1}(y(\tau))\Big)\right),
\end{gathered}
\end{equation}
where the functions $x(\tau)$ and $y(\tau)$ are defined by the equation~\eqref{eq:xy-infty1} and the functions
$\tilde{x}(\tau)$ and $\tilde{y}(\tau)$ by the equation~\eqref{eq:tide-x:tilde-y}, the functions
$\overset{\;\infty}{\tilde{A}}_{k}(\tilde{x}(\tau))$ are symmetric to the functions
$\overset{\;\infty}A_{k}(x(\tau))$ (in the sense of Remark~\ref{rem:symmetry-infty}, and the term
$(x(\tau))^k\overset{\;\infty}B_{k}(y(\tau))$ can be substituted by
$(\tilde{x}(\tau))^k\overset{\;\infty}{\tilde{B}}_{k}(\tilde{y}(\tau))$, since they are equal. Note that
if we send $N\to\infty$, then equation~\eqref{eq:u-trigAsympt-symmetric} will coincide with
equation~\eqref{eq:u-Ainfty}, because the $A$ and $B$-expansions are conjugate.

The asymptotic series~\eqref{eq:u-trigAsympt-symmetric} is convenient to use in computer calculations when
the value of the parameter $\varkappa$ is a priori unknown and is obtained as a result of computations.
\subsection{Asymptotics of the Truncated Solutions}\label{subsec:truncatedsolutions}
All Painlevé equations with an irregular singular point have so-called truncated solutions; this is a class of
special solutions whose asymptotics are given by trans-series in certain open sectors at the point at infinity.
This trans-series are obtained in the following way. Consider the asymptotic
series~\eqref{eq:app:u-asympt-regular-complete} for ``regular'' general solutions. It can be truncated in
three different ways: (1) $a_{1,1}=a_{1,-1}=0$; (2) $a_{1,1}=0$ and $a_{1,-1}\neq0$; and (3) $a_{1,-1}=0$ and
$a_{1,1}\neq0$. Such truncation means that the series~\eqref{eq:app:u-asympt-regular-complete} in case (1) does not
contain terms with the exponents $\me^{\pm\mi\theta(\tau)}$, in case (2) are abscent the terms with with
$\me^{\mi\theta(\tau)}$ and in case (3) there are no terms proportional to $\me^{-\mi\theta(\tau)}$.
Note, that in all three cases the parameter $\varkappa=0$, so that it is a characteristic property of the truncated
series/solutions.

Consider, now truncated series more carefully, the constructions of formal expansions given in
Subsections~\ref{subsec:Ainfty:trig} and \ref{subsec:trigInftyB} are valid for the truncated expansions.
Let's consider case (1). In this case we have to put in all formulae, $a_{1,\pm1}\Theta=b_{1,\pm1}=0$ and
$\varkappa=\mi b_{1,1}b_{1,-1}/3=0$ (cf. the equation~\eqref{eq:varkappab11b1m1}). Expansions
$\overset{\;\infty}{A}(x,y)$ and $\overset{\;\infty}{B}(x,y)$ coincide in this case. In each coefficient of these
expansions after we take into account conditions $b_{1,\pm1}=0$ survives only one term, namely in the
generating functions $\stackrel{\;\infty}{A}_k(x)$ it is the term proportional to $x^k$, in the generating
function $\stackrel{\;\infty}{B}_k(y)$ - the term proportional to $y^k$. Since these terms enter the corresponding
expansions in a combinations $y^k\overset{\;\infty}{A}_k(x)$ and $x^k\overset{\;\infty}{B}_k(y)$ and
$x(\tau)y(\tau)=\Theta^{-1}\tau^{-1/3}$, we should get for both expansions the same asymptotic expansion. Of course,
it is quite cumbersome to find the coefficients of the truncated asymptotic expansion by making corresponding simplifications
of the expansions for general solutions. Much easier to get this expansion from scratch using the algorithm
presented in Subsection~\ref{subsec:trigInftyB}. Note that in this case all generating functions
$\stackrel{\;\infty}{B}_{2n-1}(y)\equiv0$.
The result can be presented as the following asymptotic expnasion
\begin{equation}\label{eq:doublytruncasymptotics}
\overset{\;\infty}{A}_{trunc}(x(\tau),y(\tau))=\overset{\;\infty}{B}_{trunc}(x(\tau),y(\tau))=
1+\sum_{n=1}^\infty b_{2n}\Theta^{-2n}\tau^{-2n/3}.
\end{equation}
The first 10 coefficients of this asymptotic series are:
\begin{equation}\label{eqs:b2ndoublytruc}
\begin{gathered}
b_2=-\frac{\mi\alpha}{2},\quad
b_4=0,\quad
b_6=-\frac{\mi\alpha(\alpha^2-12)}{24},\quad
b_8=\frac{\alpha^2(\alpha^2 - 12)}{48},\\
b_{10}=\frac{3\mi\alpha(\alpha^2 - 12)}{8},\quad
b_{12}=\frac{\alpha^2(\alpha^2 - 12)(\alpha^2-66)}{144},\\
b_{14}=\frac{5\mi\alpha(\alpha^2-12)(\alpha^4-12\alpha^2-2160)}{1152},\;\;
b_{16}=-\frac{3\alpha^2(\alpha^2-12)(2\alpha^2-99)}{16},\\
b_{18}=\frac{7\mi\alpha(\alpha^2-12)(11\alpha^6-2208\alpha^4+24912\alpha^2+2721600)}{41472},\\
b_{20}=-\frac{\alpha^2(\alpha^2-12)(13\alpha^6-150\alpha^4-294264\alpha^2+13055904)}{10368}.
\end{gathered}
\end{equation}
Combination of the Wasow theorem (see Chapter 33 of \cite{Wasow1987}) with the isomonodromy results
\cite{KitVar2004} allows one to prove that under the additional assumption $\tau^{1/3}>0$ for $\tau>0$,
there exist a unique solution of the equation~\eqref{eq:dp3u} with the asymptotic
expansion~\eqref{eq:doublytruncasymptotics} in the sector $|\arg(\tau^{2/3})|<\pi$. We denote this solution
$u_0(\tau)$. Using this solution we define solutions $u_k(\tau):=-u_0(\tau\me^{\pi\mi k})$,
$k\in\mathbb{Z}$.\footnote{\label{foot:u-symmetry} Trnasformation $u(\tau)\to-u(-\tau)$ is the symmetry of the
equation~\eqref{eq:dp3u}.}
Asymptotics of these solutions are given by \eqref{eq:u-Ainfty} where
$\overset{\;\infty}{A}(x(\tau),y(\tau))=\overset{\;\infty}{B}_{trunc}(x(\tau\me^{\pi\mi k}),y(\tau\me^{\pi\mi k}))$
defined by the equation~\eqref{eq:doublytruncasymptotics} (in which $\tau\to\tau\me^{\pi\mi k}$) in the sector
$-\pi-2\pi k/3<\arg(\tau^{2/3})<\pi-2\pi k/3$. For example, on the positive real semiaxis ($\arg\,\tau=0$) there
exists exactly three solutions: $u_k(\tau)$, $k=0,\pm1$ with the truncated expansion
for the case (1).\footnote{\label{foot:defsolutionsdp3u} In \cite{KitVar2004} and all follow-up
papers by these authors, solutions to the equation~\eqref{eq:dp3u} are defined on the real positive semi-axis
($\arg\,\tau=0$) with subsequent analytic continuation. In this article, we retain this agreement.}
For $a=\mi n$ ($\alpha=2\sqrt{3} n$), $n\in\mathbb{Z}$, the functions $u_k(\tau)=u_{k+3}(\tau)$ are rational
functions of $\tau^{1/3}$, so that the series~\eqref{eq:doublytruncasymptotics} for this (and only these) choices
of the parameter $\alpha$ are convergent. For $a=0,\pm\mi$ ($\alpha=0,\mp2\sqrt{3}$), these solutions can be found
immediately from the series~\eqref{eq:doublytruncasymptotics}, because the coefficients $B_{2n}$ $|n|\geq2$
are proportional to $\alpha(\alpha^2-12)$. For $a=\mi n$, with $|n|\in\mathbb{N}$ explicit formulae for these
solutions can be found by successive application of $|n|$ B\"acklund transformations to the solution for $a=\alpha=0$.
Solutions rational in $\tau^{1/3}$ together with the B\"acklund transformations were discovered by
Gromak~\cite{Gromak1973}. Here, we note that it is important to distinguish different branches of
$\tau^{1/3}$, not only for the solutions rational in $\tau^{1/3}$, but also for the general truncated solutions,
because these solutions have different initial (and monodromy) data (see \cite{KitVar2023}, Corollary 4.1). Without
these solutions one cannot complete a description of asymptotic behaviour of solutions at the point at infinity.
We call these solutions {\it doubly truncated solutions}.

Truncations of the types (2) and (3) are mapped to each other by the symmetry considered in
Remark~\ref{rem:symmetry-infty}; this symmetry for general expansion maps the expansion to itself, in case of
the truncated solutions the symmetry is a mapping between two different truncated expansions.

The asymptotic series obtained by truncations of types (2) and (3) of the general asymptotic
series~\eqref{eq:app:u-asympt-regular-complete} are mapped onto each other by the symmetry discussed in
Remark~\ref{rem:symmetry-infty}; recall that for a general series, this symmetry maps it onto itself; in the case
of truncated series, the symmetry is a mapping between two different truncated series.
Thanks to this fact, we can translate any formula/result obtained for one of the truncated series into the
corresponding formula/result for the other. Therefore, we will consider in detail truncation of type (2).

The first thing to notice is that truncation condition of type (2) implies that $a_{k,j}=0$, $k,j\in\mathbb{N}$.
This condition means that the corresponding series depends on one parameter $a_{1,-1}$ and
does not contain the terms proportional to $\me^{\mi\theta(\tau)}$.
Thus, if we impose the condition $\mathrm{Im}\,\theta(\tau)\leq c$, with some $c\geq0$, then the series is
asymptotic. We recall that the general asymptotic series is considered in the fat ray
$|\mathrm{Im}\,\theta(\tau)|\leq c$. The exponent $\me^{-\mi\theta(\tau)}$, is decaying in the sectors of the complex
$\tau$-plane defined by the condition $\mathrm{Im}\,\theta(\tau)<0$. This condition
defines the sectors in the complex $\tau$-plane, $\Omega_0^{-}:=-\pi<\arg(\tau^{2/3})<0$.
Sector $\Omega_{-}$ is bounded by the (fat) rays where the exponent is oscillating $\me^{-\mi\theta(\tau)}$.
If we further want to apply Wasow's theory (Sections 34 and 35 of \cite{Wasow1987}), we have rearrange the original
asymptotic series~\eqref{eq:app:u-asympt-regular-complete} into the telescopic series reordering it as the
exponential series,
\begin{equation}\label{eq:telescope-0}
u^{-}_0(\tau)=\sum_{j=0}^{\infty}v_{j}^{-}(\tau)\me^{-j\mi\theta(\tau)},
\end{equation}
where in Wasow's theory $v_{0}^{-}(\tau)=u_0(\tau)$, where $u_0(\tau)$ is the solution introduced above, when
we consider truncation of the type (1); and the functions $v_{j}^{-}(\tau)$ are characterized as solutions of
the infinite system of linear ODEs. On the practical level, the functions $v_{j}^{-}(\tau)$ are obtained successively
as asymptotic series in powers of $\tau^{-1/3}$, i.e., we can construct the series for the function $v_{j}^{-}(\tau)$
only after all functions $v_{k}^{-}(\tau)$ with the subscripts $k<j$ are obtained. Following this methodology, we
get a representation of the solutions where the exponent, $\me^{-j\mi\theta(\tau)}$ appears at the first time at
the position $\omega\cdot j+1$, where $\omega$ is the standard notation for the transfinite ordinal corresponding
to the set of natural numbers. Such an asymptotic representation of truncated solutions is called a {\it transseries}.
The series for the functions
$v_{j}^{-}(\tau)$ are Borel-summable, so that these functions can be treated as the inverse Borel transform of the
corresponding Borel sums, which constitutes another approach to the transseries. We now can use
the symmetry~\footref{foot:u-symmetry} and define the series $u^{-}_k(\tau)=-u^{-}_0(\me^{\mi k}\tau)$, which
has asymptotic expansions,
\begin{equation}\label{eq:telescope-k}
u^{-}_k(\tau)=-\sum_{j=0}^{\infty}v_{j}^{-}(\me^{\mi\pi k}\tau)\me^{-\mi j\theta(\me^{\mi\pi k}\tau)},\qquad
k\in\mathbb{Z},
\end{equation}
in the sector
\begin{equation}\label{eq:Omega-m-k}
\Omega_{k}^{-}:=-\frac{3}{2}\pi-\pi k<\arg\,\tau<-\pi k.
\end{equation}
To fix the notation, we note that the truncation of type (3) of the series~\eqref{eq:app:u-asympt-regular-complete},
lead to the transseries
\begin{equation}\label{eq:telescope+k}
u^{+}_k(\tau)=-\sum_{j=0}^{\infty}v_{j}^{+}(\me^{\mi\pi k}\tau)\me^{\mi j\theta(\me^{\mi\pi k}\tau)},\qquad
k\in\mathbb{Z},
\end{equation}
in the sector
\begin{equation}\label{eq:Omega-p-k}
\Omega_{k}^{+}:=-\pi k<\arg\,\tau<-\pi k+\frac{3}{2}\pi.
\end{equation}
Here we will not dwell on further details of the construction, since they were already discussed in the
work~\cite{LinDaiTibbeoel2014}, where the authors proved that in each sector $\Omega_{k}^{\pm}$ there exists
a unique solution of the equation~\eqref{eq:dp3u} with the corresponding asymptotics; from now on, we will denote
(with some abuse in notation) these solutions in the same way as the corresponding transseries, $u^{\pm}_k(\tau)$.

As an example of the application of this classical theory, we see that on the positive real semi-axis ($\arg\,\tau=0$) there exist
exactly two solutions whose asymptotics are given by the transseries, $u^{-}_{-1}(\tau)$ and $u^{+}_{+1}(\tau)$;
these solutions (and the corresponding transseries) depend on one parameter $a_{1,-1}$ or $a_{1,1}$, respectively.
In \cite{Var2025} parameters were expressed in terms of the monodromy data.

A specific feature of the truncated solutions to the Painlev\'e equations is that one can present them in terms
of the truncated asymptotic series~\eqref{eq:app:u-asympt-regular-complete} which approximate these solutions not
only in the open sectors $\Omega_{k}^{\pm}$, but also in their closure, $\bar{\Omega}_{k}^{\pm}$, in the (fat) rays
representing the boundaries of the sectors, the truncated solutions are oscillating, rather than exponentially decay.
The parameter $a_{1,-1}$ or, respectively, $a_{1,1}$ characterizing a truncated solution is the same in
$\bar{\Omega}_{k}^{\pm}$. Here, we will not discuss/use this feature and concentrate on $A/B$-expansions for these
solutions.

Applying the scheme of Subsection~\ref{subsec:trigInftyB}, I calculated first 13 polynomials
$\overset{\;\infty}{B}_k(y)$ in the truncated case (2), $a_{1,1}=\Theta^{-1}b_{1,1}=0$. The polynomials
in the truncated case are considerably simpler than those for the general solutions, so that this calculation
takes about 10 minutes. The polynomials are
even functions of $y$ with $\deg\,\overset{\;\infty}{B}_{2n}(y)=\deg\,\overset{\;\infty}{B}_{2n+1}(y)=2n$. The first
eight polynomials read:
\begin{align}
\overset{\;\infty}{B}_0(y)=&1,\quad
\overset{\;\infty}{B}_1(y)=b_{1,-1},\quad
\overset{\;\infty}{B}_2(y)=-\frac{\mi\alpha y^2}{2}+\frac{b_{1,-1}^2}{3},\label{eqs:Btrig1-7}\\
\overset{\;\infty}{B}_3(y)=&-\frac{\mi(\alpha^2+4\alpha-1)b_{1, -1}y^2}{8}+\frac{b_{1,-1}^3}{12},\nonumber\\
\overset{\;\infty}{B}_4(y)=&-\frac{\mi(\alpha^2+4\alpha+1/3)b_{1,-1}^2y^2}{12}+\frac{b_{1,-1}^4}{54},\nonumber\\
\overset{\;\infty}{B}_5(y)=&-\frac{(\alpha+3)^2(\alpha-1/3)(\alpha-3)b_{1,-1}y^4}{2^7}
-\frac{\mi(\alpha^2+4\alpha+1/9)b_{1,-1}^3y^2}{2^5}\nonumber\\
&+\frac{5b_{1,-1}^5}{6^4},\nonumber\\
\overset{\;\infty}{B}_6(y)=&-\frac{\mi\alpha(\alpha^2-12)y^6}{24}
-\frac{(3\alpha^4+16\alpha^3-10\alpha^2-128\alpha-57)b_{1,-1}^2y^4}{2^53^2}\nonumber\\
&-\frac{\mi(\alpha^2+4\alpha-1/6)b_{1,-1}^4y^2}{2^23^3}+\frac{b_{1,-1}^6}{6^4},\nonumber\\
\overset{\;\infty}{B}_7(y)=
&\mi\big(\alpha^6-4\alpha^5-67\alpha^4-152\alpha^3+579\alpha^2+2076\alpha-225\big)\frac{b_{1,-1}y^6}{2^{10}3}\nonumber\\
&-\Big(9\alpha^4+56\alpha^3+10\alpha^2-328\alpha-\frac{329}{3}\Big)\frac{b_{1,-1}^3y^4}{2^93}\nonumber\\
&-\mi\Big(25\alpha^2+100\alpha-\frac{29}{3}\Big)\frac{b_{1,-1}^5y^2}{2^73^4}+\frac{7b_{1,-1}^7}{6^6}.\nonumber
\end{align}
The expansion~\eqref{eq:AxyB} for $B_k(y)=\overset{\;\infty}{B}_k(y)$ and $x=x(\tau)$ and $y=y(\tau)$
defined by the equations~\eqref{eq:infty:regular-w-theta-tau} and \eqref{eq:xy-infty1} (with $\nu+1=\varkappa=0$
for the truncated solutions):
$x(\tau)y(\tau)=\Theta^{-1}\tau^{-1/3}$ and $x(\tau)=\Theta^{-1}\tau^{-1/3}\me^{-\mi\theta(\tau)}$, exactly
coincide with the truncated asymptotic series~\eqref{eq:app:u-asympt-regular-complete}.

Let us now consider the $A$-expansion~\eqref{eq:A-inftyDef}. At first glance, the function
$y(\tau)=\me^{\mi\theta(\tau)}$ in $\Omega_0^{-}$ increases, so the expansion~\eqref{eq:A-inftyDef} has only a
formal meaning. However, the functions $\overset{\;\infty}{A}_k(x)$ in the truncated case $\varkappa=0$ are
proportional to $x^k$, therefore the $k$-th term of the expansion, $y^k\overset{\;\infty}{A}_k(x)$, after
the substitution $x=x(\tau)$, $y=y(\tau)$ has the order $\tau^{-k/3}$ as $\tau\to\infty$, $\tau\in\Omega_0^{-}$,
and the expansion~\eqref{eq:A-inftyDef} is in fact an asymptotic series. Thus, to form the asymptotic series
in the truncated case it is convenient to rewrite the expansion~\eqref{eq:A-inftyDef} in the following way,
\begin{equation}\label{eq:A-infty-truc-Def}
\overset{\;\infty}{A}_{trunc}(x,y)=
\sum_{k=0}^{\infty}(yx)^{2k}\left.\frac{\overset{\;\infty}{A}_{2k}(x)}{x^{2k}}\right\vert_{\varkappa=0}.
\end{equation}
Recall that the functions $\overset{\;\infty}{A}_{2k+1}(x)=0$, as proved in Remark~\ref{rem:Bn-even}.
Below, we will focus on constructing the asymptotic series. The reader interested in explicit
formulae for the coefficients of the truncated $B$-expansion should substitute $\varkappa=0$ into the corresponding
formulae for the coefficients of the general solution; see equations~\eqref{eq:ann} (does not contain $\varkappa$, cf.
$\overset{\;\infty}{B}_n(0)$, $n=0,\ldots,6$ in \eqref{eqs:Btrig1-7}), \eqref{eq:an-(n-2)}, and \eqref{eq:an4-n}.

Now let us consider explicit formulas for the first few terms of the expansion~\eqref{eq:A-infty-truc-Def}.
The first term of the expansion~\eqref{eq:A-infty-truc-Def} is the function $\overset{\;\infty}{A}_{0}(x)$,
which does not depend on $\varkappa$ and is therefore given by the same equation~\eqref{eq:A0-inftyRat} as in
the general case. However, the next term undergoes some simplifications,
\begin{equation}\label{eq:A2-trunk-parfrac}
\begin{gathered}
\left.\frac{\overset{\;\infty}{A}_{2}(x)}{x^{2}}\right\vert_{\varkappa=0}=
-\frac{\mi b_{1,-1}x}{12}-\frac{\mi\alpha}{2}-6\mi-\frac{\mi(3\alpha^2+12\alpha-5)}{2(1-b_{1,-1}x/6)^3} \\
+\frac{3\mi(3\alpha^2+12\alpha-13)}{4(1-b_{1,-1}x/6)^2}-\frac{\mi(3\alpha^2+12\alpha-53)}{4(1-b_{1,-1}x/6)},
\end{gathered}
\end{equation}
\begin{equation}\label{eq:A4-trunk-parfrac}
\begin{gathered}
\left.\frac{\overset{\;\infty}{A}_{4}(x)}{x^{4}}\right\vert_{\varkappa=0}=-\frac{b_{1,-1}x}{48}\Bigg(
\frac{b_{1,-1}x}{9}+\frac{3\alpha^2+36\alpha+107}{6}\\
-\frac{9\alpha^4+88\alpha^3+234\alpha^2+152\alpha+165}{4(1-xb_{1,-1}/6)^3}
+\frac{9\alpha^4+72\alpha^3+114\alpha^2-120\alpha+25}{4(1-xb_{1,-1}/6)^4}\\
+\frac{9\alpha^4+120\alpha^3+618\alpha^2+1272\alpha+493}{24(1-xb_{1,-1}/6)^2}\Bigg).
\end{gathered}
\end{equation}
\begin{equation}\label{eq:A6-trunk-parfrac}
\begin{gathered}
\left.\frac{\overset{\;\infty}{A}_{6}(x)}{x^{6}}\right\vert_{\varkappa=0}=
\frac{\mi b_{1,-1}^3x^3}{2^83^4}+\frac{\mi(3\alpha^2+36\alpha+85)b_{1,-1}^2x^2}{2^63^4}\\
+\frac{\mi(3\alpha^4+56\alpha^3+238\alpha^2-328\alpha-2289)b_{1,-1}x}{2^93^2}\\
-\frac{\mi(\alpha^3+18\alpha^2+60\alpha-30)}{24}
+\frac{\mi}{192}\sum_{k=1}^5\frac{\mu_k}{(1-b_{1,-1}x/6)^k},
\end{gathered}
\end{equation}
where
\begin{equation*}\label{eq:mu1}
\mu_1=\frac{1}{2^3}(3\alpha^6+84\alpha^5+1041\alpha^4+7240\alpha^3+31141\alpha^2+63204\alpha+15383),
\end{equation*}
\begin{equation*}\label{eq:mu2}
\mu_2=-\frac{1}{2^3}(45\alpha^6+876\alpha^5+7167\alpha^4+31288\alpha^3+81483\alpha^2+116828\alpha+\frac{111787}{3}),
\end{equation*}
\begin{equation*}\label{eq:mu3}
\mu_3=\frac{1}{2^2}(75\alpha^6+1188\alpha^5+7257\alpha^4+21576\alpha^3+33865\alpha^2+28996\alpha+\frac{24689}{3}),
\end{equation*}
\begin{equation*}\label{eq:mu4}
\mu_4=-\frac{1}{2}(45\alpha^6+612\alpha^5+2871\alpha^4+5208\alpha^3+2769\alpha^2+540\alpha-\frac{2875}{3}),
\end{equation*}
\begin{equation*}\label{eq:mu5}
\mu_5=9\alpha^6+108\alpha^5+387\alpha^4+216\alpha^3-645\alpha^2+300\alpha-\frac{125}{3}.
\end{equation*}
To calculate $\left.\frac{\overset{\;\infty}{A}_{6}(x)}{x^{6}}\right\vert_{\varkappa=0}$ using the scheme in
Subsection~\ref{subsec:Ainfty:trig}, it is necessary to calculate the leading coefficient of the polynomial
$\overset{\;\infty}{B}_7(y)$ presented for this purpose above.

Combining together equations~\eqref{eq:u-Ainfty}, \eqref{eq:A-infty-truc-Def}, and \eqref{eq:xy-infty1}, we obtain
the following asymptotic series for the truncated solutions in the sector $\Omega_0^{-}$,
\begin{equation}\label{eq:TruncsolAsympserOmega-m}
\begin{gathered}
u(\tau)\underset{\genfrac{}{}{0pt}{}{\tau\to\infty}{\tau\in\Omega_0^{-}}}{=}
\frac{\varepsilon(\varepsilon b)^{2/3}}{2}\tau^{1/3}
\sum_{k=0}^{\infty}\Theta^{-2k}\tau^{-2k/3}
\left.\frac{\overset{\;\infty}{A}_{2k}(x)}{x^{2k}}\right\vert_{\varkappa=0},\\
x=\Theta^{-1}\tau^{-1/3}\me^{-\mi\Theta^2\tau^{2/3}},
\end{gathered}
\end{equation}
where the scaling parameter $\Theta=3^{3/4}(\varepsilon b)^{1/6}>0$. The asymptotic
series~\eqref{eq:TruncsolAsympserOmega-m} is significantly more complicated than the original truncated
asymptotic series~\eqref{eq:app:u-asympt-regular-complete}; further research is needed to understand
whether this asymptotic series provides any new practical advantages, such as a better approximation of $u(\tau)$
near the boundaries of the $\Omega_0^{-}$ sector.
\section{Irregular Singularity: Elliptic Expansions}\label{sec:elliptic}
The explicit form of the complete "elliptic" asymptotic expansions of the Painlev\'e equations are not known yet.
That makes a difficulty in confirmation of the expansions~\eqref{eq:uAxy} and \eqref{eq:AxyB} in this case.
Recently an extensive studies of the Boutroux-type elliptic asymptotics for the third-fifth Painlev\'e equations
were undertaken by
S. Shimomura~\cite{Shimomura2022P5,Shimomura2024EllErr,Shimomura2024DP3,Shimomura2025}.
These results concerns the
leading term of the corresponding asymptotics as well as the estimate for the correction terms. The latter results
should be compared with the formal constructions presented in this section.

To simplify our constructions in this section, we now assume that $a=0$.
We begin with the definition of the domain $\Omega$, in which we will construct our formal expansions.
Put $\tau=r\me^{\mi\varphi_0}$, where $\varphi_0\neq\pi k/2$, $k\in\mathbb{Z}$, is some fixed angle in the
complex $\tau$-plane and $r\in\mathbb{C}$. On the other hand, we can present $\tau$ via the polar coordinates,
$\tau=|\tau|\me^{\mi\varphi}$, so that $\mathrm{Re}(r^{2/3})=|\tau|^{2/3}\cos(2(\varphi-\varphi_0)/3)$,
$\mathrm{Im}(r^{2/3})=|\tau|^{2/3}\sin(2(\varphi-\varphi_0)/3)$. Now, define the domain $\Omega$ as
a "fat" ray originating from the point at infinity at an angle $-\varphi_0$, more precisely,
\begin{equation}\label{def:Omega-elliptic}
\Omega:=\{\tau:\;|\mathrm{Im}(r^{2/3})|<\delta_0, |r|>\eta_0,\;\mathrm{where}\;\delta_0>0\;\mathrm{and}\;\eta_0>0\}.
\end{equation}
In this section, we will discuss two types of function pairs $(x(\tau), y(\tau))$ that are related to elliptic Jacobi
functions with modulus $\kappa$.\footnote{\label{foot:Jacobi} We use notation for these functions adopted in
\cite{BE3}.}
\subsection{Elliptic Expansion: Non-conjugate Variables}\label{subsec:e-non-conjugate}
The first straightforward generalization of the functions $x(\tau)$, $y(\tau)$ introduced in the trigonometric
case (see equations ~\eqref{eq:infty:regular-w-theta-tau} and ~\eqref{eq:xy-infty1}), reads
\begin{equation}\label{eq:x-y-elliptic1}
\begin{aligned}
x(\tau)&=\Theta^{-1}r^{-1/3}\operatorname{dn}(\vartheta(r)/2,\kappa)^{-1/2}
\me^{-2\mi\operatorname{am}(\vartheta(r)/2,\kappa)},\\
y(\tau)&=\me^{2\mi\operatorname{am}(\vartheta(r)/2,\kappa)},
\end{aligned}
\end{equation}
where $\Theta$ is defined in \eqref{eq:xy-infty1} and
\begin{equation}\label{eq:vartheta}
\vartheta(r):=\Theta^2 Pr^{2/3}+\vartheta_0,\qquad
P,\vartheta_0\in\mathbb{C}.
\end{equation}
The points $r$ where $\operatorname{dn}(\vartheta(r)/2,\kappa)=0$ should be excluded from this construction,
we, however, do not pay any special attention to this points, because one can make a similar construction with
$\operatorname{dn}(\vartheta(r)/2,\kappa)^{+1/2}$, so that zeros of $\operatorname{dn}(\vartheta(r)/2,\kappa)=0$
does not play any role in the theory of the function $u(\tau)$.

Define differential operator $D$ acting in the space of analytic functions of two variables $x$ and $y$
\begin{equation}\label{eq:D-elliptic-def}
D:=\frac{\mi P}{x^2y}\frac{\partial}{\partial y}
-\left(\frac{x}{2}+ \frac{\mi P}{xy^2}\left(1 + \frac{\kappa^2(y^2 - 1)}{4(4y+\kappa^2(y - 1)^2)}\right)\right)
\frac{\partial}{\partial x}.
\end{equation}
Consider the following partial differential equation
\begin{equation}\label{eq:A-inftyPDE-elliptic}
D^2\ln\,\overset{\;\infty}{A}_e^1=-\frac{p^2}{3x^4y^4\left(1+\kappa^2(y - 1)^2/(4y)\right)}
\left(\overset{\;\infty}{A}_e^1
-\frac{1}{\Big(\overset{\;\infty}{A}_e^1\Big)^2}\right),
\end{equation}
where $D$ is defined in \eqref{eq:D-elliptic-def}, $p=\me^{2\mi\varphi_0/3}$, and
$\overset{\;\infty}{A}_e^1\equiv\overset{\;\infty}{A}_e^1(x,y)$.\footnote{\label{foot:Ainfty1ek}
The superscript $1$ is used in the notation of the function $\overset{\;\infty}{A}^1_{e}(x)$,
since in Subsection~\ref{subsec:e-conjugate}, we will consider the generating function
$\overset{\;\infty}{A}_e^2(x,y)$, which is related to a different way of introducing elliptic variables, moreover,
the construction similar to that used for the function $\overset{\;\infty}{A}_e^2(x,y)$ also works for the variables
$x(\tau)$, $y(\tau)$ considered in this subsection.}
If a solution $u(\tau)$ of the equation~\eqref{eq:dp3u} admits representation~\eqref{eq:uAxy} with rational functions
$A_k(x)$, more precisely, in the elliptic case we employ the following notation
\begin{equation}\label{eq:A-infty-elliptic-def}
u(\tau)=\frac{\varepsilon(\varepsilon b)^{2/3}}{2}\tau^{1/3}\overset{\;\infty}{A}_e^1(x(\tau),y(\tau)),
\end{equation}
where the functions $x(\tau)$ and $y(\tau)$ are defined in the equations~\eqref{eq:x-y-elliptic1}, and
\begin{equation}\label{eq:Ainfty-expansion1}
\overset{\;\infty}{A}_e^1(x,y)=\sum_{k=0}^{\infty} y^k\overset{\;\infty}{A}^1_{e,k}(x),
\end{equation}
where $\overset{\;\infty}{A}^1_{e,k}(x)$, $k=0,1,2,\ldots$, are rational functions
of $x$.
Our immediate goal is to demonstrate the existence of a formal solution~\eqref{eq:Ainfty-expansion1}
to PDE~\eqref{eq:A-inftyPDE-elliptic}.

Consider construction of the coefficients $\overset{\;\infty}{A}^1_{e,k}(x)$. Substituting the
expansion~\eqref{eq:Ainfty-expansion1} into the equation~\eqref{eq:A-inftyPDE-elliptic}, expanding the resulting
expression in the power-series in $y$ and successively equating
to $0$ the corresponding coefficients, we get obtain ODEs defining the coefficients
$\overset{\;\infty}{A}^1_{e,k}(x)$. One finds the first equation defining $\overset{\;\infty}{A}^1_{e,0}(x)$, by
equating to $0$ the coefficient of the term $y^{-4}$; this is a nonlinear ODE, which after simple transformations
reads
\begin{equation}\label{eq:ODE-Ainfty-elliptic-coeff0}
3x\overset{\;\infty}{A}^1_{e,0}(x)\frac{\md^2}{\md x^2}\overset{\;\infty}{A}^1_{e,0}(x)
-3x\left(\frac{\md}{\md x}\overset{\;\infty}{A}^1_{e,0}(x)\right)^2
+5\overset{\;\infty}{A}^1_{e,0}(x)\frac{\md}{\md x}\overset{\;\infty}{A}^1_{e,0}(x)=0.
\end{equation}
The general solution of the equation~\eqref{eq:ODE-Ainfty-elliptic-coeff0} is
\begin{equation}\label{eq:Ainfty-elliptic-coeff0gen}
\overset{\;\infty}{A}^1_{e,0}(x)= q\exp(-3c/(2x^{2/3})),\qquad
q\in\mathbb{C}\setminus\{0\},\quad
c\in\mathbb{C}.
\end{equation}
Thus the only rational solution of \eqref{eq:ODE-Ainfty-elliptic-coeff0} is
\begin{equation}\label{eq:Ainfty-elliptic-coeff0q}
\overset{\;\infty}{A}^1_{e,0}(x)= q.
\end{equation}
The ODEs for the rest coefficients $\overset{\;\infty}{A}^1_{e,k}(x)$, $k\in\mathbb{N}$, are linear inhomgeneous
second order ODEs, whose inhomogeneous part is a polynomial in $x^2$; more precisely,
\begin{equation}\label{eq:LODE-Ainfty-elliptic-coeffs-k}
\begin{gathered}
P^{2k}\kappa^{2k}q^{2k-1}\!\left(\!9x^2\frac{\md^2}{\md x^2}\overset{\;\infty}{A}^1_{e,k}(x)
-3(8k-5)x\frac{\md}{\md x}\overset{\;\infty}{A}^1_{e,k}(x)\right.\\
\left.+8k(2k-1)\overset{\;\infty}{A}^1_{e,k}(x)\!\right)=R_k(x,P,p,q,\kappa),
\end{gathered}
\end{equation}
where the polynomial $R_k\in\mathbb{R}[x,P,p,q,\kappa]$ is an even function of $x$ and
$\deg_x\,R_k=2\lfloor\frac{2(k-1)}{3}\rfloor$.\footnote{\label{foot:def-floor} The notation $\lfloor r \rfloor$
represents the integer part of the rational number $r$.}

Let's consider the first few equations~\eqref{eq:LODE-Ainfty-elliptic-coeffs-k}:\\
(1) $k=1$,
\begin{equation}\label{eq:LODE-A1infty-elliptic}
P^2\kappa^2q\left(9x^2\frac{\md^2}{\md x^2}\overset{\;\infty}{A}^1_{e,1}(x)
-9x\frac{\md}{\md x}\overset{\;\infty}{A}^1_{e,1}(x)+8\overset{\;\infty}{A}^1_{e,1}(x)\right)=\frac{64}{3}p^2(q^3-1).
\end{equation}
Equation~\eqref{eq:LODE-A1infty-elliptic} has the only one particular rational solution, namely,
\begin{equation}\label{eq:A1infty-elliptic}
\overset{\;\infty}{A}^1_{e,1}(x)=\frac{8p^2(q^3-1)}{3P^2\kappa^2q};
\end{equation}
(2) $k=2$,
\begin{equation}\label{eq:LODE-A2infty-elliptic}
P^4\kappa^4q^3\left(9x^2\frac{\md^2}{\md x^2}\overset{\;\infty}{A}^1_{e,2}(x)
-33x\frac{\md}{\md x}\overset{\;\infty}{A}^1_{e,2}(x)+48\overset{\;\infty}{A}^1_{e,2}(x)\right)=.R_2,
\end{equation}
where
\begin{equation*}\label{eq:R2}
R_2=\frac{256}{9}p^2(q^3-1)(3P^2q^2(\kappa^2-2)+2p^2(4q^3-1)),
\end{equation*}
The general rational solution of the equation~\eqref{eq:LODE-A2infty-elliptic} reads
\begin{equation}\label{eq:A2infty-elliptic}
\overset{\;\infty}{A}^1_{e,2}(x)=q_2x^2+R_2/(48P^4\kappa^4q^3),
\quad
q_2\in\mathbb{C};
\end{equation}
(3) $k=3$,
\begin{equation}\label{eq:LODE-A3infty-elliptic}
P^6\kappa^6q^5\!\left(\!9x^2\frac{\md^2}{\md x^2}\overset{\;\infty}{A}^1_{e,3}(x)
-57x\frac{\md}{\md x}\overset{\;\infty}{A}^1_{e,3}(x)
+120\overset{\;\infty}{A}^1_{e,3}(x)\!\right)\!=R_3,
\end{equation}
where
\begin{equation}\label{eq:R3}
\begin{gathered}
R_3=R_{3,2}x^2+R_{3,0},\\
R_{3,2}=\frac{16}{3}P^3q^3\kappa^4\big(3P^3q^2q_2(2-\kappa^2)+4Pp^2q_2(4q^3-1)+8\mi p^2q(q^3-1)\big),\\
R_{3,0}=\frac{64}{81}p^2(q^3-1)\big(9P^4q^4(23\kappa^4-128\kappa^2+128)\\
+240P^2p^2q^2(4q^3\kappa^2-8q^3-\kappa^2+2)+32p^4(34q^6-26q^3+1)\big).
\end{gathered}
\end{equation}
The general rational solution of equation~\eqref{eq:LODE-A3infty-elliptic} can be presented in the following form
\begin{equation}\label{eq:A3infty-elliptic}
\overset{\;\infty}{A}^1_{e,3}(x)=q_4x^4+\frac{R_{3,2}x^2}{24P^6\kappa^6q^5}+\frac{R_{3,0}}{120P^6\kappa^6q^5},
\quad
q_4\in\mathbb{C};
\end{equation}
We see that for $k=1,2$, and $3$, $\overset{\;\infty}{A}^1_{e,k}(x)$ are, in fact, are polynomials of $x^2$. Since
for the functions $R_k$, one can get a recurrence relation in terms of $\overset{\;\infty}{A}^1_{e,m}(x)$ with
$m<k$, one proves the following statement via the mathematical induction argument:
\begin{equation}\label{eq:Rk}
R_k=\sum_{m=0}^N R_{k,2m}x^{2m},\qquad
R_{k,2m}=R_{k,2m}(P,p,q,\kappa),\qquad
N=\left\lfloor\frac{2(k-1)}{3}\right\rfloor,
\end{equation}
\begin{equation}\label{eq:Ak:Dkm}
\begin{gathered}
\overset{\;\infty}{A}^1_{e,k}(x)=q_{2(N+1)}x^{2(N+1)}+\frac{1}{P^{2k}\kappa^{2k}q^{2k-1}}
\sum_{m=0}^N\frac{R_{k,2m}}{D_k(m)}x^{2m},\\
D_k(m)=4(2k-3m)(2k-3m-1),
\end{gathered}
\end{equation}
where $q_{2(N+1)}=0$, if $D_k(N+1)\neq0$ and $q_{2(N+1)}\in\mathbb{C}$, if $D_k(N+1)=0$. By other words,
$q_{2(N+1)}=0$ for $k=1+3l$, $l\in\mathbb{Z}_{\geqslant0}$, otherwise $q_{2(N+1)}\in\mathbb{C}$.

Thus, we have constructed a general rational formal solution of the partial differential
equation~\eqref{eq:A-inftyPDE-elliptic}. To obtain the corresponding solution of the equation~\eqref{eq:dp3u},
we need to choose the parameters $P$, $q$, $\kappa$ and $q_{2n}$, $n\in\mathbb{N}$.
Below, in the equation~\eqref{eq:Ainfty-e-Ox2}, we operate with a formal series, and the notation
$\mathcal{O}(x^2)$ means all terms proportional to $x^2$:\footnote{\label{foot:xto0}Note that for $\tau\to\infty$
$(x(\tau))^2=\mathcal{O}(\tau^{-2/3})$.}
\begin{align}
\overset{\;\infty}{A}_e^1(x,y)&=
\sum_{k=0}^{\infty} y^k\big(\overset{\;\infty}{A}^1_{e,k}(0)+\mathcal{O}(x^2)\big)=
\sum_{k=0}^{\infty} y^k\overset{\;\infty}{A}^1_{e,k}(0)+\mathcal{O}(x^2)\nonumber\\
=q+\overset{\;\infty}{A}^1_{e,1}(0)
&\sum_{k=1}^{\infty} y^k\frac{\overset{\;\infty}{A}^1_{e,k}(0)}{\overset{\;\infty}{A}^1_{e,1}(0)}+\mathcal{O}(x^2)
=q+\frac{8p^2(q^3-1)}{3P^2\kappa^2q}\sum_{k=1}^{\infty} ky^k+\mathcal{O}(x^2)\label{eq:Ainfty-e-Ox2}\\
&=q+\frac{8p^2(q^3-1)}{3P^2\kappa^2q}\frac{y}{(1-y)^2}+\mathcal{O}(x^2),\label{eq:Ainfty-e-main}
\end{align}
where in the equation~\eqref{eq:Ainfty-e-Ox2} we replaced $\overset{\;\infty}{A}^1_{e,1}(0)$ with its
expression~\eqref{eq:A1infty-elliptic} and set
\begin{equation}\label{eq:Ainftyk:Ainfty1}
\overset{\;\infty}{A}^1_{e,k}(0)/\overset{\;\infty}{A}^1_{e,1}(0)=k\quad
\mathrm{for}\quad
k\in\mathbb{N}.
\end{equation}
The infinite system of equations~\eqref{eq:Ainftyk:Ainfty1} impose the following conditions on the parameters:
\begin{equation}\label{eq:kappaPp-qs}
\kappa^2=\frac{(s-1)(s+3)}{(s+1)(s-3)},\quad
\frac{P^2}{p^2}=\frac{2q}{3}\frac{(s-3)}{(s-1)},\qquad
s^2=8q^3+1.
\end{equation}
Substituting into the equation~\eqref{eq:Ainfty-e-main}: (1) relations~\eqref{eq:kappaPp-qs}; and (2) for $x$ and $y$,
variables $x(\tau)$ and $y(\tau)$, respectively, which are defined by the equations~\eqref{eq:x-y-elliptic1};
using the definition of the elliptic sine, after the straightforward simplifications, one, recalling the definition
of the function $u(\tau)$ \eqref{eq:A-infty-elliptic-def}, arrives at the following
``candidate'' for a leading term of the large-$\tau$ asymptotics:

Now, we obtain ``our candidate'' for the leading term of the large-$\tau$ asymptotics of the function $u(\tau)$.
For this purpose, using the relations~\eqref{eq:kappaPp-qs}, we simplify the formula~\eqref{eq:Ainfty-e-main},
\begin{equation}\label{eq:Ainfty-e-main-simple}
\overset{\;\infty}{A}_e^1(x,y)=q+4\frac{(s-3)}{(s-1)}\frac{qy}{(1-y)^2}+\mathcal{O}(x^2).
\end{equation}
Substituting in the equation~\eqref{eq:Ainfty-e-Ox2} for $x$ and $y$ the variables $x(\tau)$ and $y(\tau)$,
respectively, which are defined by equations~\eqref{eq:x-y-elliptic1} and using a definition of the elliptic sine,
we derive for $u(\tau)$ (see equation~\eqref{eq:A-infty-elliptic-def}) the following asymptotics,
\begin{equation}\label{eq:u-asympt-tau-infty-e}
u(\tau)\underset{\tau\to\infty,\arg\tau=\varphi_0}{=}\frac{\varepsilon(\varepsilon b)^{2/3}}{2}\tau^{1/3}
\left(q-\frac{(s-3)}{(s-1)}\frac{q}{\operatorname{sn}^2(\vartheta(r)/2,\kappa)}+\mathcal{O}\big(\tau^{-2/3}\big)\right),
\end{equation}
where $\vartheta(r)$ is defined by the equation~\eqref{eq:vartheta} and $\tau=r\me^{\mi\varphi_0}$ with
$r\in\Omega$ (see equation~\eqref{def:Omega-elliptic}).
Moreover, after formula~\eqref{eq:u-asympt-tau-infty-e} is obtained, the fat ray $\Omega$ should be turned into
a ``cheese-shaped'' domain by cutting out disks around the zeros of the elliptic sine. Depending on the radius of
these disks, the error estimate in formula~\eqref{eq:u-asympt-tau-infty-e} may deteriorate. We will not discuss
these standard asymptotic issues here, since our main goal is to consider the algebraic aspects of deriving the
formula for the leading term of the asymptotics of $u(\tau)$.

The formula~\eqref{eq:u-asympt-tau-infty-e} contains two parameters $P$ and $\vartheta_0$, because due to
the relations~\eqref{eq:kappaPp-qs}, $\kappa$ is given in terms of $P$ and $p=\me^{2\mi\varphi_0/3}$.
At first glance, this may not seem bad, as we are dealing with the second-order ODE \eqref{eq:dp3u},
however, in fact, this means that not all the parameters have been determined, because elliptic asymptotics
of the Painlev\'e equations along a fixed ray in the complex plane contain only one parameter in their leading
terms.\footnote{\label{foot:Boutroux}
This was first shown by P. Boutroux~\cite{Boutroux}, who consider the example of the first Painlev\'e equation.}
The second parameter appears only when we consider $\vartheta_0$ as a function of $\varphi_0$.
In principle, we haven't studied yet the $y$-series that appear as the coefficients of the powers of $x^{2n}$,
the coefficients of latter series depend on $P$ and the other parameters $q_2, q_4,\ldots$, so that their
convergence can impose some additional conditions on $P$. One of the instruments to study this question is
the conjugate expansion~\eqref{eq:AxyB}, however, in this case the functions $x(\tau)$ and $y(\tau)$ are not
conjugate, at least for generic values of $P$, so that the structure of the expansion~\eqref{eq:AxyB} is not clear
yet.

On the other hand, the Boutroux-type asymptotics for the degenerate third Painlev\'e equation was recently obtained by
S. Shimomura~\cite{Shimomura2024DP3} via the isomonodromy deformation method. In the next subsection,
we will examine in detail the correspondence between the asymptotics obtained by Shimomura and the asymptotic
formula~\eqref{eq:u-asympt-tau-infty-e}, that allows us to fix the parameter $P$ in terms of the so-called
modulus function defined by Shimomura.
\subsection{Comparison with the Boutroux-Type Asymptotics}\label{subsec:shimomura}
In this subsection, we use the subscript $s$ in the notation $x$ and $y$ introduced in Shimomura's paper
\cite{Shimomura2024DP3}, because the notation $x$ and $y$ used in the previous section have a
different meaning. Shimomura considered the degenerate third Painlev\'e equation in the following form,
\begin{equation}\label{eq:dp3s}
y_s''=\frac{(y_s')^2}{y_s}-\frac{y_s'}{x_s}-2y_s^2+\frac{3a}{x_s}+\frac{1}{y_s},
\end{equation}
where $y_s=y_s(x)$ and $y_s'=\md y_s(x_s)/\md x_s)$, which comes from \eqref{eq:dp3u} via the substitution
 \begin{equation}\label{eq:s-u-variables}
 \varepsilon\tau u(\tau)= (x_s/3)^2y_s(x_s),\qquad
 \varepsilon b\tau^2 = 2(x_s/3)^3.
\end{equation}
Recall that throughout this section, we assume, for simplicity, that $a=0$, the latter assumption having no effect
on the leading term of the asymptotics of $y_s(x_s)$, unless we consider a parametrization of this asymptotics
by monodromy data and corresponding connection formulae for the asymptotics, which we do not do here.

In this subsection, we assume that the reader is familiar with the standard notation and well-known facts from
the theory of the Weierstrass elliptic $\wp$-function; all the necessary information is collected in \cite{BE3} and
will be used hereafter without references.
The asymptotic formula obtained in \cite{Shimomura2024DP3} reads
\begin{equation}\label{eq:dp3s-e-asympt}
y_s(x_s)=\wp\big(\mi(x_s-x_0^+)+\mathcal{O}(x_s^{-\delta});g_2,g_3\big)+\frac{A_\phi}{12},
\end{equation}
where $x_s=t_s\me^{\mi\phi}\to\infty$, such that $0<\phi<\pi/3$, and $t_s\in\mathbb{C}$ belongs to a ``cheese-like''
strip defined in \cite{Shimomura2024DP3}; $\delta>0$, the phase-shift $x_0$ is given explicitly via monodromy data
of the associated linear ODE, and
\begin{equation}\label{eq:g2g3}
g_2=\frac{(A_{\phi})^2}{12},\qquad
g_3=\frac{(A_{\phi})^3}{6^3}-1.
\end{equation}
The complex-valued modulus function $A_{\phi}$ is defined in a certain way via a normalization condition for
an elliptic integral. The reader is referred to the paper~\cite{Shimomura2024DP3} for precise definitions.

Our main goal in this subsection is to establish the correspondence between the asymptotic
formulae~\eqref{eq:dp3s-e-asympt} and \eqref{eq:Ainfty-e-main-simple}.
Equations~\eqref{eq:A-infty-elliptic-def} and \eqref{eq:s-u-variables} imply that
\begin{equation}\label{eq:A-y}
\overset{\;\infty}{A}_e^1(x(\tau),y(\tau))=2^{1/3}y_s(x_s).
\end{equation}
Let us rewrite asymptotics of the function $\overset{\;\infty}{A}_e^1(x(\tau),y(\tau))$ in terms of the
Weierstrass $\wp$-function,
\begin{equation}\label{eq:Ainfty-e-W}
\begin{aligned}
\overset{\;\infty}{A}_e^1(x(\tau),y(\tau))&=
q-\frac{(s-3)}{(s-1)}\frac{q}{\mathrm{sn}^2(\vartheta(r)/2,\kappa)}+\mathcal{O}\big(\tau^{-2/3}\big)\\
&=q-q\frac{(s-3)}{(s-1)}\frac{\wp(\mi(\tilde{x}-\tilde{x}_0))-e_3}{e_1-e_3}+\mathcal{O}\big(\tau^{-2/3}\big),
\end{aligned}
\end{equation}
where there exists $\mu\in\mathbb{C}\setminus\{0\}$, such that
\begin{equation}\label{eq:e1e2e3}
e_1=(2-\kappa^2)\mu^2,\qquad
e_2=(2\kappa^2-1)\mu^2,\quad
e_3=-(\kappa^2+1)\mu^2,
\end{equation}
and
\begin{equation}\label{eq:vartheta-tilde-x}
\vartheta(r)/2=(e_1-e_3)^{1/2}\mi(\tilde{x}-\tilde{x_0})=\sqrt{3}\,\mi\mu(\tilde{x}-\tilde{x_0}).
\end{equation}
The equation~\eqref{eq:vartheta-tilde-x} together with the definition of $\vartheta(r)$
(see equation~\eqref{eq:vartheta}) implies
that
\begin{equation}\label{eq:tilde-x-x0}
\tilde{x}=-\frac{\mi\Theta^2P}{2\sqrt{3}\mu}r^{2/3},\qquad
\vartheta_0=-2\sqrt{3}\,\mi\mu\tilde{x}_0.
\end{equation}
Using the second equation~\eqref{eq:s-u-variables}, we find
\begin{equation}\label{eq:xs-tau}
x_s=\frac{3}{2^{1/3}}(\varepsilon b)^{1/3}\tau^{2/3}=\frac{\Theta^2}{2^{1/3}\sqrt{3}}r^{2/3}p
=t_s\me^{\mi\phi},\quad
p=\me^{2\mi\varphi_0/3}\Rightarrow
\phi=2\varphi_0/3,
\end{equation}
where in the second equation the definition of $\Theta$ (see equation~\eqref{eq:xy-infty1}) is used.
The compatibility condition for the Weierstrass $\wp$-functions in \eqref{eq:dp3s-e-asympt} and \eqref{eq:Ainfty-e-W}
reads
\begin{equation}\label{eq:compatibility-e}
x_s=\tilde{x}\quad
\mathrm{and}\quad
\tilde{x}_0=x_0^+;
\end{equation}
as follows from the
equations~\eqref{eq:tilde-x-x0} and \eqref{eq:xs-tau}, it is equivalent to the following relations,
\begin{equation}\label{eq:mu-vartheta0}
\vartheta_0=-2\sqrt{3}\,\mi\mu x_0^+,\qquad
\mu=-\frac{\mi}{2^{2/3}}\frac{P}{p}.
\end{equation}
Now, the second equation~\eqref{eq:kappaPp-qs} allows to calculate $\mu^2$ in terms of the elliptic coordinates,
\begin{equation}\label{eq:mu-s-q}
\mu^2=-\frac{q}{2^{1/3}3}\frac{(s-3)}{(s-1)}.
\end{equation}
We proceed now with the simplification of the leading term of the asymptotics~\eqref{eq:Ainfty-e-W}; it can be
rewritten as follows,
\begin{gather}\label{eq:leading-term-Ainfty-e-W}
q+q\frac{(s-3)}{(s-1)}\frac{e_3}{e_1-e_3}-q\frac{(s-3)}{(s-1)}\frac{\wp(\mi(\tilde{x}-\tilde{x}_0))}{e_1-e_3}=
q\left(1-\frac{(s-3)}{(s-1)}\frac{(\kappa^2+1)}{3}\right)\\
-\frac{q(s-3)}{3\mu^2(s-1)}\wp(\mi(\tilde{x}-\tilde{x}_0))=\frac{q}{3}\frac{(s^2+3)}{(s^2-1)}
+2^{1/3}\wp(\mi(x_s-x_0^+)),
\end{gather}
where we used the equations~\eqref{eq:e1e2e3}, \eqref{eq:kappaPp-qs}, \eqref{eq:mu-s-q}, and
\eqref{eq:compatibility-e}. Substituting now for the functions in both sides of the equation~\eqref{eq:A-y}, their
asymptotics, we arrive at the relation between the modulus function and the parameter $q$,
\begin{equation}\label{eq:Aphi-q}
2^{1/3}\frac{A_{\phi}}{12}=\frac{q}{3}\frac{(s^2+3)}{(s^2-1)}\quad
\Rightarrow\quad
A_{\phi}=2^{2/3}\frac{2q^3+1}{q^2}.
\end{equation}
Denoting $q=2^{1/3}z$, one finds that $4z^3-A_{\phi}z^2+1=0$, so that parameter $2^{1/3}q$ has a sense of a branching
point of the elliptic curve $w^2=4z^3-A_{\phi}z^2+1=0$, which is used by Shimomura for the definition of the
function $A_{\phi}$. Since the invariants $g_2$ and $g_3$ are well-known symmetric functions of the roots
$e_1$, $e_2$, and $e_3$ (cf. \cite{BE3}, the relations~\eqref{eq:mu-s-q}, \eqref{eq:e1e2e3}, and
\eqref{eq:Aphi-q}, allow one to reproduce equations~\eqref{eq:g2g3}.

The second equation~\eqref{eq:Aphi-q} solves our problem, it fixes the parameters: $q$, $s$ and therefore $P$ and
$\kappa$ for a given $\varphi_0$ (see the third equation~\eqref{eq:xs-tau}) and allows one to use
equation~\eqref{eq:u-asympt-tau-infty-e} as the asymptotic formula for $u(\tau)$; the error estimate in this formula
is written under the assumption that the parameter $q_2, q_4,\ldots $ can be chosen such that the corresponding
series in $y$ converges in the same cheese-shaped region as that defined for the leading term of the asymptotics.
Since these series are not studied yet, the reader should rely on the less sharp error estimate obtained by
Shimomura \eqref{eq:dp3s-e-asympt}.
\subsection{Elliptic Expansion: Symmetry-Conjugate Variables}\label{subsec:e-conjugate}
To define the conjugate variables $\hat{x}$ and $\hat{y}$ in the elliptic case, we use the same notation as in
Subsection~\ref{subsec:e-non-conjugate}, but change the definition of the variables $x$ and $y$ to a more symmetric
form,
\begin{equation}\label{eq:hatx-haty-e}
\begin{aligned}
x(\tau)&
=\Theta^{-1}r^{-1/3}\mathrm{dn}^{-1/2}(\vartheta(r)/2,\kappa)\me^{-\mi\,\mathrm{am}(\vartheta(r)/2,\kappa)},\\
y(\tau)&
=\Theta^{-1}r^{-1/3}\mathrm{dn}^{-1/2}(\vartheta(r)/2,\kappa)\me^{\mi\,\mathrm{am}(\vartheta(r)/2,\kappa)}.
\end{aligned}
\end{equation}
For this choice of the functions $x(\tau)$ and $y(\tau)$ the differential operator
\begin{equation}\label{eq:D-e-conjugate}
\begin{aligned}
D&=
\left(\frac{y}{2}+\frac{\mi P}{x}\left(-\frac{1}{2}+\frac{\kappa^2(y^2-x^2)}{4(4xy+\kappa^2(y - x)^2)}\right)\!
\right)\frac{\partial}{\partial y}\\
&+\left(\frac{x}{2}+\frac{\mi P}{y}\left(\frac{1}{2}+\frac{\kappa^2(y^2-x^2)}{4(4xy+\kappa^2(y - x)^2)}\right)\!
\right)\frac{\partial}{\partial x},
\end{aligned}
\end{equation}
and the PDE defining $A/B$-expansions is
\begin{equation}\label{eq:A-inftyPDE-e-conjugate}
D^2\ln\,\overset{\;\infty}{A}_e^2=-\frac{4p^2}{3xy(4xy+\kappa^2(y-x)^2)}
\left(\overset{\;\infty}{A}_e^2
-\frac{1}{\Big(\overset{\;\infty}{A}_e^2\Big)^2}\right),
\end{equation}
where $\overset{\;\infty}{A}_e^2=\overset{\;\infty}{A}_e^2(x,y)$ or $\overset{\;\infty}{B}_e^2(x,y)$
(see below).
For the functions $\overset{\;\infty}{A}_e^2(x,y)$ and $\overset{\;\infty}{B}_e^2(x,y)$, one can find two types of
conjugate formal expansions: one of them has the standard form~\eqref{eq:uAxy}, \eqref{eq:AxyB}, which is considered
in this section and the other to be discussed elsewhere.
Our goal is to construct the most general formal solutions of the equation~\eqref{eq:A-inftyPDE-e-conjugate} that
are defined by the following expansions:
\begin{equation}\label{eq:Ainfty-e2-regular-A-B-expansions}
\overset{\;\infty}{A}_e^2(x,y)=\sum_{k=0}^{\infty} y^k\overset{\;\infty}{A}^2_{e,k}(x),\qquad
\overset{\;\infty}{B}_e^2(x,y)=\sum_{k=0}^{\infty} x^k\overset{\;\infty}{B}^2_{e,k}(y),
\end{equation}
where $\overset{\;\infty}{A}^2_{e,k}(x)$, and $\overset{\;\infty}{B}^2_{e,k}(y)$ are rational functions of $x$ and
$y$, respectively. As series of two independent variables $x$ and $y$, the formal solutions
$\overset{\;\infty}{A}_e^2(x,y)$ and $\overset{\;\infty}{B}_e^2(x,y)$ are different. However, they are related to
each other by symmetry; more precisely, the following statement holds.
\begin{proposition}\label{prop:ABP}
Assume that the expansion of $\overset{\;\infty}{A}^2_{e}(x)$ given by the first
equation~\eqref{eq:Ainfty-e2-regular-A-B-expansions}, with coefficients
$\overset{\;\infty}{A}^2_{e,k}(x)\equiv\overset{\;\infty}{A}^2_{e,k}(x,P)$ that are rational functions of $x$, is
a formal solution of the equation~\eqref{eq:A-inftyPDE-e-conjugate}.
Define the functions
\begin{equation}\label{eq:ABP}
\overset{\;\infty}{B}^2_{e,k}(y,P):=\overset{\;\infty}{A}^2_{e,k}(y,-P),\qquad
k\in\mathbb{Z}_{\geqslant0}.
\end{equation}
Then the expansion of $\overset{\;\infty}{B}_e^2(x,y)$ given by the second
equation~\eqref{eq:Ainfty-e2-regular-A-B-expansions}
with coefficients $\overset{\;\infty}{B}^2_{e,k}(y)=\overset{\;\infty}{B}^2_{e,k}(y,P)$, which are rational
functions of $y$, is a formal solution of the equation~\eqref{eq:A-inftyPDE-e-conjugate},
and vice versa.
\end{proposition}
Let's consider the expansion $\overset{\;\infty}{B}_e^2(x,y)$ and find the first few rational functions:
\begin{equation}\label{eq:B2e0B2e1}
\overset{\;\infty}{B}^2_{e,0}(y)=q,\quad
q\in\mathbb{C}\setminus\{0\},\qquad
\overset{\;\infty}{B}^2_{e,1}(y)=c_2y^3+c_1y+\frac{8p^2(q^3-1)}{3\kappa^2P^2qy},
\end{equation}
\begin{equation}\label{eq:B2e2}
\begin{gathered}
\overset{\;\infty}{B}^2_{e,2}(y)=c_4y^6+c_3y^4\\
+\frac{3\kappa^2P^2qc_1^2+6\mi\kappa^2Pq^2c_1-12(\kappa^2-2)P^2q^2c_2+16p^2(2q^3+1)c_2}{6\kappa^2P^2q^2}\,y^2
\\
-\frac{2\big(3(\kappa^2-2)P^3q^2c_1-4Pp^2(4q^3-1)c_1+8\mi p^2q(q^3-1)\big)}{9\kappa^2P^3q^2}\\
+\frac{16p^2(q^3 - 1)\big(3(\kappa^2-2)P^2q^2+2p^2(4q^3-1)\big)}{27\kappa^4P^4q^3y^2}.
\end{gathered}
\end{equation}
The following proposition can be proved with the help of mathematical induction argument based on the
equations~\eqref{eq:B2e0B2e1} and \eqref{eq:B2e2}.
\begin{proposition}\label{prop:B2e}
Let the partial differential equation~\eqref{eq:A-inftyPDE-e-conjugate} have parameters
$P,p,\kappa\in\mathbb{C}\setminus\{0\}$. Assume that $\overset{\;\infty}{A}_e^2$, given by the second
equation~\eqref{eq:Ainfty-e2-regular-A-B-expansions}, is its formal solution
with $\overset{\;\infty}{B}^2_{e,0}(y)\not\equiv0$ and rational (in $y$) coefficient functions
$\overset{\;\infty}{B}^2_{e,k}(y)$. Then
\begin{equation}\label{eq:Binfty2-k-m}
\overset{\;\infty}{B}^2_{e,0}(y)=q\in\mathbb{C}\setminus\{0\}
\quad\mathrm{and}\quad
\overset{\;\infty}{B}^2_{e,k}(y)=\sum_{m=0}^{2k}\overset{\;\infty}{B}^{2,2m-k}_{e,k} y^{2m-k},\quad
k\in\mathbb{N},
\end{equation}
where the coefficients $\overset{\;\infty}{B}^{2,2m-k}_{e,k}$ are independent of $y$
rational functions of $P$, $p^2$, $\kappa^2$ and for $k\geqslant1$ are polynomials in
$c_1,c_2,\ldots,c_{2k-1},c_{2k}\in\mathbb{C}$.
Moreover,\footnote{\label{foot:prod} Below we assume that $\prod_{1}^0=1$.}
\begin{gather}
\overset{\;\infty}{B}^{2,3k}_{e,k}=c_{2k}y^{3k},\quad
\overset{\;\infty}{B}^{2,3k-2}_{e,k}=c_{2k-1}y^{3k-2},\quad
c_{2k},c_{2k-1}\in\mathbb{C},\qquad
k\in\mathbb{N},
\label{eq:Binfty-e-k-3k+0-2}\\
\overset{\;\infty}{B}^{2,2m-k}_{e,k}=
\sum_{\genfrac{}{}{0pt}{}{\sum jn_j\leqslant m}{n_j\geqslant0}}
\mathfrak{b}_{k,n_1,\ldots,n_{2k-2}}(\kappa^2,P,p,q)\prod_{j=1}^{2k-2}c_j^{n_j},
\;\;
m=0,1,\ldots,2k-2,
\end{gather}
where the coefficients $\mathfrak{b}_{k;n_1,\ldots,n_{2k-2}}(\kappa^2,P,p,q)$ are rational functions of their
variables. In particular, for $k\geqslant2$,
\begin{equation}\label{eq:B2ekmk:B2ek2mk}
\begin{aligned}
\overset{\;\infty}{B}^{2,-k}_{e,k}&=\mathfrak{b}_{k;0,\ldots,0}(\kappa^2,P,p,q),\\
\overset{\;\infty}{B}^{2,2-k}_{e,k}&=\mathfrak{b}_{k;1,\ldots,0}(\kappa^2,P,p,q)c_1+
\mathfrak{b}_{k;0,\ldots,0}(\kappa^2,P,p,q),
\end{aligned}
\end{equation}
Conversly, let $p,P,\kappa\in\mathbb{C}\setminus\{0\}$ and $c_{n}\in\mathbb{C}$, $n\in\mathbb{N}$, then there
exists a formal solution $\overset{\;\infty}{A}_e^2$ given by the second
equation~\eqref{eq:Ainfty-e2-regular-A-B-expansions} with $\overset{\;\infty}{B}^2_{e,-1}(y)\not\equiv0$ whose
coefficients $\overset{\;\infty}{B}^2_{e,k}(y)$ are rational functions of $y$, $P$, $p^2$, and $\kappa^2$.
These coefficients have representation~\eqref{eq:Binfty2-k-m}. Under the normalization~\eqref{eq:Binfty-e-k-3k+0-2}
the solution is unique.
\end{proposition}
\begin{corollary}\label{prop:A2e}
Under the same assumptions as in Proposition~{\rm\ref{prop:B2e}} there exists the unique formal solution
$\overset{\;\infty}{A}^2_{e}$ of
the partial differential equation~\eqref{eq:A-inftyPDE-e-conjugate} given by the first
equation~\eqref{eq:Ainfty-e2-regular-A-B-expansions}, with $\overset{\;\infty}{A}^2_{e,-1}(x)\not\equiv0$ and rational
in $x$, $P$, $p^2$, $\kappa^2$ coefficient functions $\overset{\;\infty}{A}^2_{e,k}(x)$. Denoting explicitly the
dependence on $P$ of the coefficients $\overset{\;\infty}{A}^2_{e,k}(x)=\overset{\;\infty}{A}^2_{e,k}(x;P)$ and
$\overset{\;\infty}{B}^2_{e,k}(y)=\overset{\;\infty}{B}^2_{e,k}(y;P)$, one have the following relation
$\overset{\;\infty}{A}^2_{e,k}(x;P)=\overset{\;\infty}{B}^2_{e,k}(x;-P)$.
\end{corollary}

Now, consider a relation of the expansion $\overset{\;\infty}{B}_e^2(x,y)$ to the construction of the Boutroux-type
asymptotics of the solution $u(\tau)$. In this case, the relation~\eqref{eq:A-infty-elliptic-def}, should be
substituted by
\begin{equation}\label{eq:B-infty-elliptic-def}
u(\tau)=\frac{\varepsilon(\varepsilon b)^{2/3}}{2}\tau^{1/3}\overset{\;\infty}{B}_e^2(x(\tau),y(\tau)),
\end{equation}
where the functions $x(\tau)$ and $y(\tau)$ are defined by the equations~\eqref{eq:hatx-haty-e}. This construction
is similar to the one that was used in Subsection~\ref{subsec:e-non-conjugate} for the function
$\overset{\;\infty}{A}_e^1(x,y)$; we sum all terms in the formal solution $\overset{\;\infty}{B}_e^2(x,y)$,
which, after substituting $x(\tau)$ and $y(\tau)$ for the variables $x$ and $y$, respectively, have the behavior
$\mathcal{O}(1)$ as $\tau\to\infty$;
\begin{equation}\label{eq:B2e-main}
\begin{aligned}
\overset{\;\infty}{B}_e^2(x,y)&=q+
\sum_{k=1}^{\infty}\left(\!\!\left(\frac{x}{y}\right)^k\overset{\;\infty}{B}^{2,-k}_{e,k}
+\mathcal{O}\left(\frac{x^k}{y^{k-2}}\right)\!\!\right)\\
&=q+\overset{\;\infty}{B}^{2,-1}_{e,1}
\sum_{k=1}^{\infty} \left(\frac{x}{y}\right)^k
\frac{\overset{\;\infty}{B}^{2,-k}_{e,k}}{\overset{\;\infty}{B}^{2,-1}_{e,1}}+\mathcal{O}(x^2)\\
&=q+\frac{8p^2(q^3-1)}{3\kappa^2P^2q}\sum_{k=1}^{\infty} k\left(\frac{x}{y}\right)^k
+\mathcal{O}(x^2)\\
&=q+\frac{8p^2(q^3-1)}{3\kappa^2P^2q}\frac{x/y}{(1-x/y)^2}+\mathcal{O}(x^2),
\end{aligned}
\end{equation}
where in the last equation~\eqref{eq:B2e-main} we assume that the parameters $P$ and $\kappa$ are chosen such that
they solve the following system of equations,
\begin{equation}\label{eq:B2ek:B2e1}
\overset{\;\infty}{B}^{2,-k}_{e,k}/\overset{\;\infty}{B}^{2,-1}_{e,1}=k\quad
\mathrm{for}\quad
k\in\mathbb{N}.
\end{equation}
Since $\overset{\;\infty}{B}^{2,-k}_{e,k}=\overset{\;\infty}{A}^{1}_{e,k}(0)$, the solution of this system
coincide with the solution of the system~\eqref{eq:Ainftyk:Ainfty1} and given by the same
equations~\eqref{eq:kappaPp-qs}. {\bf In the rest of this subsection, we assume that the parameters
$P$, $p$, $\kappa^2$, and $q$ satisfy the conditions given in equations~\eqref{eq:kappaPp-qs}.}
Now, we can reproduce the leading term of the asymptotics
\begin{equation}\label{eq:B2ek-asymp}
\begin{aligned}
\overset{\;\infty}{B}_e^2(x(\tau),y(\tau))&=
q+\frac{8p^2(q^3-1)}{3\kappa^2P^2q\big(\sqrt{y(\tau)/x(\tau)}-\sqrt{x(\tau)/y(\tau)}\big)^2}
+\mathcal{O}\big(\tau^{-2/3}\big)\\
&=q-\frac{(s-3)}{(s-1)}\frac{q}{\operatorname{sn^2(\vartheta(r)/2,\kappa)}}
+\mathcal{O}\big(\tau^{-2/3}\big),
\end{aligned}
\end{equation}
where we used the conditions~\eqref{eq:kappaPp-qs} and the definitions~\eqref{eq:hatx-haty-e} of the functions
$x(\tau)$ and $y(\tau)$. Now, combining together equations ~\eqref{eq:B-infty-elliptic-def} and
\eqref{eq:B2ek-asymp}, one arrives again at the asymptotics~\eqref{eq:u-asympt-tau-infty-e} for the function $u(\tau)$.

Consider the correction term in the equation~\eqref{eq:B2e-main},
\begin{equation}\label{eq:x^2-corr}
\mathcal{O}\big(x^{2}\big)=\sum_{k=1}^{\infty}\overset{\;\infty}{B}^{2,2-k}_{e,k}\frac{x^k}{y^{k-2}}
+\mathcal{O}\big(x^{4}\big),
\end{equation}
where the notation $\mathcal{O}\big(x^{4}\big)$ represents all terms of the second
expansion \eqref{eq:Ainfty-e2-regular-A-B-expansions}, which, when $y=x$ is substituted, have the order $x^4$.
Substituting now into the equation~\eqref{eq:x^2-corr} $x(\tau)$ and $y(\tau)$, determined by the
equations~\eqref{eq:hatx-haty-e}, we obtain the correction term in the equation~\eqref{eq:B2ek-asymp}
in the form of a one-sided complex Fourier series with respect to the amplitude function,
\begin{equation}\label{eq:corrterm-Fourier}
\begin{gathered}
\mathcal{O}\big(\tau^{-2/3}\big)=\Theta^{-2}\tau^{-2/3}
\Big(\operatorname{dn}\big(\vartheta(r)/2,\kappa\big)\Big)^{-1}
\sum_{k=1}^{\infty}b_{e,k}^{2-k}\me^{-2i(k-1)\operatorname{am}\big(\vartheta(r)/2,\kappa\big)}\\
+\mathcal{O}\big(\tau^{-4/3}\big),
\end{gathered}
\end{equation}
where
\begin{equation}\label{eq:def-b-ek-2mk}
\tau^{2/3}=pr^{2/3},\qquad
b_{e,k}^{2-k}=
\left.p\overset{\;\infty}{B}^{2,2-k}_{e,k}\right\vert_{eq.\eqref{eq:kappaPp-qs}}.
\end{equation}
The last notation means that the coefficients $\overset{\;\infty}{B}^{2,2-k}_{e,k}$ are calculated for the parameters
satisfying the relations~\eqref{eq:kappaPp-qs}. According to the last equation~\eqref{eq:B2ekmk:B2ek2mk},
the coefficients $b_{e,k}^{2-k}$ are linear functions of $c_1$. The coefficients can be presented as follows:
\begin{equation}\label{eq:be11:be20}
b_{e,1}^{1}=c_1p,\quad
b_{e,2}^{0}=-4\mi\frac{P}{p}+\frac{(3s^2-2s-9)}{(s+3)(s-1)}\,2c_1p,
\end{equation}
where
\begin{equation*}
\frac{P}{p}=\frac{(s^2-1)^{1/6}}{\sqrt{3}}\frac{\sqrt{s-3}}{\sqrt{s-1}},
\end{equation*}
and for $k\geqslant2$
\begin{equation}\label{eq:bek2-k}
\begin{gathered}
b_{e,k}^{2-k}=-4\mi\frac{P}{p}\frac{Q_{2k-4}(s)}{(2k-1)!(s+3)^{k-2}(s-1)^{k-2}}\\
+\frac{(2(k-1)^2+1)(s-3)(s+1)+4s}{(s+3)(s-1)}\,2c_1p,
\end{gathered}
\end{equation}
where $Q_{2k-4}(s)$ are polynomials in $s$ with $\deg(Q_{2k-4}(s))=2k-4$ and integer coefficients. The first six
polynomials $Q_{2k-4}(s)$ read:
\begin{equation*}\label{eq:polynomialsQ-elliptic}
\begin{aligned}
k=2:&\;\;Q_0(s)=6,\quad
k=3:\;\;Q_2(s)=2^63^2(s^2-5),\quad\\
k=4:&\;\;
Q_4(s)=2^43(1241s^4+2436s^3-13266s^2-11340s+16065),\\
k=5:&\;\;
Q_6(s)=2^{12}3^3(71s^6+284s^5-585s^4-3264s^3+2037s^2\\
&+3780s-2835),\\
k=6:&\;\;
Q_8(s)=2^83^3(197823s^8+1184920s^7-824620s^6-11736120s^5\\
&-16680534s^4+47405160s^3+2899764s^2-46153800s+22920975),\\
k=7:&\;\;
Q_{10}(s)=2^{13}3^3(1398679s^{10}+11187384s^9+7677153s^8\\
&-147117360s^7-79588362s^6+30783168s^5+721058418s^4\\
&-875438928s^3-222469821s^2+1313512200s-492567075).
\end{aligned}
\end{equation*}
As follows from \eqref{eq:bek2-k} a part of the series~\eqref{eq:corrterm-Fourier} proportional to $c_1$ is
easy to sum up, however, the rest of the series looks complicated and most probably related with a solution of
an equation of the Lam\'e type. The latter equation, together with the value of the parameter $c_1$,
could be obtained by a substitution of the expansion~\eqref{eq:B2ek-asymp} with $\mathcal{O}\big(\tau^{-2/3}\big)$
given by equation~\eqref{eq:corrterm-Fourier} into the degenerate third Painlev\'e equation~\eqref{eq:dp3u}.

If instead of the expansion~\eqref{eq:B-infty-elliptic-def}, we consider the conjugate expansion
\begin{equation}\label{eq:A2-infty-elliptic-def}
u(\tau)=\frac{\varepsilon(\varepsilon b)^{2/3}}{2}\tau^{1/3}\overset{\;\infty}{A}_e^2(x(\tau),y(\tau)),
\end{equation}
then, because of the symmetry stated in Corollary~\ref{prop:A2e}, we arrive at: the same
constraints~\eqref{eq:kappaPp-qs} for the parameters $P$, $p$, $q$, and $s$; exactly the same formula for
the leading term of the asymptotic expansion $\overset{\;\infty}{A}_e^2(x(\tau),y(\tau))$
as for $\overset{\;\infty}{B}_e^2(x(\tau),y(\tau))$ (see the second equation~\eqref{eq:B2ek-asymp}); and
to the formula for the correction term $\mathcal{O}\big(\tau^{-2/3}\big)$, which differs from the corresponding
formula for the function $\overset{\;\infty}{B}_e^2(x(\tau),y(\tau))$ (see equation~\eqref{eq:corrterm-Fourier})
by the reflection of $P\to-P$ in $\vartheta(r)$ and in all the coefficients $b_{e,k}^{2-k}$.
The sign of the parameter
$P$ is not fixed in our construction, so we will not discuss this issue at this stage; it may be related
to the representation of the correction term in different sectors where the elliptic asymptotics is valid.

\end{document}